\documentclass{amsart}
\usepackage{amsthm, a4wide}
\usepackage{amsfonts}
\usepackage{amsmath}
\usepackage{amssymb}
\usepackage[english]{babel}
\usepackage{verbatim}
\usepackage{enumerate}

\usepackage{mathtools}

\usepackage[toc,page]{appendix}

\usepackage{color}

\let\ol=\overline

\numberwithin{equation}{section}
\newtheorem{thm}{Theorem}[section]
\newtheorem{lem}[thm]{Lemma}
\newtheorem{corol}[thm]{Corollary}
\newtheorem{prop}[thm]{Proposition}

\theoremstyle{definition} 
\newtheorem{defn}[thm]{Definition}
\newtheorem{example}[thm]{Example}

\theoremstyle{remark}
\newtheorem{rem}[thm]{Remark}

\newcommand{\R}{\mathbb{R}}

\newcommand{\N}{{\mathbb N}}

\newcommand{\Cc}{{\mathcal C}}

\newcommand{\g}{\gamma}
\newcommand{\mC}{{\mathcal C}}
\let\d=\delta
\let\e=\varepsilon 
\let\l=\lambda
\let\s=\sigma
\let\r=\rho 
\let\o=\omega
\let\ol=\overline
\let \a =\alpha
\let \t =\tau
\let \b =\beta

\let\Si=\Sigma
\newcommand{\rh}{\hat{\r}}
\newcommand{\rb}{\ol{\r}}

\newcommand{\mA}{\mathcal{A}} 
\newcommand{\mF}{\mathcal{F}}

\newcommand{\mR}{\mathcal{R}}

\newcommand{\mT}{\mathcal{T}}

\def\Diff{{\it Diff}}
\def\bel{\begin{equation}\label}
\def\eeq{\end{equation}}
\def\pT{\partial\mathcal{T}}

\allowdisplaybreaks[4]

\begin{document}

\title[Regularity for degenerate eikonal equations]{Regularity of the minimum time and of \\ viscosity solutions of degenerate eikonal equations\\ 
via generalized Lie brackets}

\author [Bardi, Feleqi, and Soravia]{Martino Bardi$^1$, \; Ermal Feleqi$^2$, \; {and} \; Pierpaolo Soravia$^3$}
 \address{1 and 3: Dipartimento di Matematica "T. Levi-Civita",
Universit\`a degli Studi di Padova,
via Trieste, 63, I-35121 Padova, Italy; 
2: Department of Mathematics, University of Vlora, Albania.}
\email{bardi@math.unipd.it,  ermal.feleqi@univlora.edu.al, soravia@math.unipd.it}

\begin{abstract}
In this paper we relax the current regularity theory for the eikonal equation by using the recent theory of { set-valued} iterated Lie brackets. We give sufficient conditions for small time local attainability of general, symmetric, nonlinear systems, which have as a consequence the H\"older regularity of the minimum time function in optimal control. We then apply such result to prove H\"older continuity of solutions of the Dirichlet boundary value problem for the eikonal equation with low regularity of the coefficients.
We also prove that the sufficient conditions for the H\"older regularity are essentially necessary, at least for smooth vector fields and target.
\end{abstract}%
\subjclass{}
\keywords{Time-optimal control, minimum time functions, nonlinear systems, Lie brackets, geometric control theory, eikonal equation, Bellman equations, Hamilton-Jacobi equations.
 {\em MR Subject Classification:} 35F30 (49L25 93B05 93B27).
}

 \thanks{
The first author is partially supported by the research project ``Nonlinear Partial Differential Equations: Asymptotic Problems and Mean-Field Games" of the Fondazione CaRiPaRo. 
He is also member of the Gruppo Nazionale per l'Analisi Matematica, la Probabilit{\`a} e le loro Applicazioni (GNAMPA) of the Istituto Nazionale di Alta Matematica (INdAM)
}
\maketitle
\section{Introduction}

In this paper we 
address the question of regularity of viscosity solutions of the Dirichlet boundary value problem for degenerate eikonal equations,  { namely, }
\begin{equation}\label{eqstart}
\left\{ \begin{aligned}
& \sum_{i=1}^m| f_i u|^2 + 2\sum_{i=1}^m b_i(x)  f_iu = h^2(x)\, ,  \quad \mbox{ in } \Omega\\  
&u =g \qquad  \mbox{ on } \partial \Omega,
\end{aligned} \right.  
\end{equation}
under appropriate compatibility of the boundary condition. Here $b_i,h$ are given 
 coefficients, and $f_i$ are a family of vector fields, written in coordinates as differential operators $f_i  = \sum_{j=1}^n f_i^j(x) \partial_{x_j}$, so that $f_i u(x) = \sum_{j=1}^n f_i^j(x) \partial_{x_j}u(x)$. The solution $u$ will be continuous up to the boundary and meant as viscosity solution. { The Hamilton-Jacobi equation in \eqref{eqstart} is degenerate at some point  $x\in \bar\Omega$ if at such point the vectors $f_i(x)$ do not span all $\R^n$ and so the Hamiltonian is not coercive in the moment variables. 
We are interested in equations with such degeneracies, especially at boundary points.}

Under appropriate regularity of the coefficients in the equation, it is known that if 
 the Hamiltonian is coercive with respect to the gradient of the solution $u$, then $u$ is locally Lipschitz continuous,  and that this is false in general when such a property is not satisfied. In this case one 
 { aims at the H\"older regularity} of the solutions with a suitable 
  exponent, and properties of the Lie algebra generated by the vector fields $f_i$ come into play. Therefore the regularity of the vector fields is a key assumption. For a review of the classical theory in this direction we refer the reader to the book \cite{BaCa97} { and the references therein}.  It is also well known that continuous viscosity solutions are unique and they have a representation formula as value functions of an appropriate control problem. In the case of the homogeneous boundary conditions ($g\equiv0$) and positive and constant Lagrangian (e.g. $h\equiv1$), then the solution is the minimum time function in optimal control from the target $\R^n\backslash\Omega$, see \cite{BaCa97,BaSo91}. It is also known that the H\"older regularity of the minimum time function is a consequence of the small time local attainability of the target by the family of vector fields, more precisely of suitable estimates of the minimum time function with powers of the distance function from $\partial\Omega$. Such estimates can be derived from properties of the iterated Lie brackets between the available vector fields, see, for instance, 
 \cite{so} and \cite{kr}. 

In a classical setting an iterated Lie bracket of length $k+1$ 
is defined for vector fields at least 
 of class $C^{k+1}$. 
 For instance, $[f_i,f_j]=Df_j\;f_i-Df_i\;f_j$ is a length one Lie bracket and we need to be able to compute continuously the Jacobians of the two vector fields $f_i,f_j$ so that it is a continuous vector field. { Rampazzo and Sussmann introduced a notion of commutator 
 for non-smooth vector fields by defining the first bracket for Lipschitz fields, and studied the properties of the flows ~\cite{RaSu01, RaSu07}. Their approach was continued by Feleqi and Rampazzo \cite{FeRa1} who defined iterated Lie brackets of length 
   $k+1$ 
   if the vector fields are of class $C^k$ and the $k$th-derivatives are locally Lipschitz continuous. In this context a Lie bracket of length $k+1$ is } 
{ defined in the classical way} { only almost everywhere and it is completed as a multivalued map elsewhere. Such theory of non-smooth commutators is growing and aims at recovering several classical results 
  of geometric control \cite{FeRa2,FeRa}.  Applications to mechanical or engineering problems have not yet been pursued, but they are expected in view of the the well-known fruitful interactions of nonsmooth analysis and control theory  (see, e.g., \cite{Vin}).}

{ The main goal of this paper 
 is to extend the current existence and regularity} theory for the eikonal equation when 
 it is necessary to use 
such multivalued iterated Lie brackets. However, some results are new also in the case of smooth vector fields.  

 We will 
first study  the regularity of the minimum time function { to reach the boundary $\partial\Omega$}, a key tool to reach regularity of solutions of (\ref{eqstart}). We will do so 
 for families of vector fields which are fully nonlinear, therefore in a 
 { wider} generality than what is needed for the mere equation (\ref{eqstart}). We also drop regularity of the boundary of the domain $\partial\Omega$ by allowing $\R^n\backslash\Omega$ to be the union of the closure of an open set and a locally finite set of isolated points, 
 whereas  around points of $\partial\Omega$ that are not isolated $\partial\Omega$ needs not to be smooth but just satisfy an exterior cone property.  { At such points we assume the existence of a generalized bracket pointing outward $\Omega$. In PDE terms this is a weak form of non-characteristicity of the boundary for the eikonal equation.}
We 
 obtain $1/m$-H\"older regularity of the minimum time function, 
where $m$ is the highest length of the brackets involved.
We then apply such result to the solution of the   Dirichlet problem \eqref{eqstart} for the eikonal equation with right-hand side $h\ne 0$ and boundary data $g$ continuous and satisfying natural compatibility conditions. { Assuming only the bracket condition at boundary points, the problem is well-posed in the subset $\mR$ of $\Omega$ of the positions from where the boundary can be reached, with the solution going to $+\infty$ at points of $\partial\mR$, and the well-posedness is in the whole $\Omega$ if the generalized brackets satisfy a nonsmooth H\"ormander condition everywhere.}
  Finally, for smooth data, we show that the sufficient condition on the Lie brackets for the $1/2$-H\"older continuity of the minimum time function becomes also necessary if completed with the possibility of exiting $\Omega$ by means of a single tangential vector field. 
    To our knowledge this kind of necessary conditions is completely new in the literature. 

Small time local attainability and regularity of the minimum time function is a long studied and important subject in optimal control.
Classical results by Petrov \cite{pe} show sufficient conditions for attainability at a single point by requiring that the convex hull of the vector fields at the point contains the origin in its interior. This is called a first order controllability condition. Liverovskii \cite{li} studied the corresponding problem of second order when a similar request is made on the family of vector fields augmented with their first 
 Lie brackets, see also Bianchini and Stefani \cite{bist}. Controllability of higher order to a point was studied by Liverovskii \cite{li2}. For attainability of a target different from a point we recall the papers by Bacciotti \cite{ba} in the case of targets of codimension 1, { by one of the authors and Falcone for piecewise smooth targets \cite{BaFa90}, and by another author of this paper} for manifolds of any dimension and possibly with a boundary \cite{so}. { The typical approach of these papers consists of showing that the signed distance function from the boundary of the target becomes negative along suitable admissible trajectories; this remains the starting point of the proofs in most of the following 
 papers, including the present one.}
 Monti and Serra Cassano \cite{MSC} used tools of geometric measure theory to prove that the Carnot-Carath\'eodory distance solves the PDE in \eqref{eqstart} almost everywhere in a suitable sense.
{ Trelat \cite{Tre} studied the sub-analiticity of the sub-Riemannian distance and of viscosity solutions of the problem \eqref{eqstart} in the case of analytic vector fields $f_i$ and subanalytic $\Omega$ and $g$.}
More recently the work by Krastanov and Quincampoix \cite{kr,kr2,kr3} pointed out the importance of the geometry of the target and studied higher order attainability with smooth families of vector fields but nonsmooth targets, for affine systems with nontrivial drift.  For the same class of systems Marigonda, Rigo and Le \cite{ma,ma2,ma3, ma4} studied higher order regularity focusing on the lack of smoothness of the target and the presence of state constraints. { A regularity result for the solution of a very special case of \eqref{eqstart} was given in \cite{Al12} and we will discuss it in Section \ref{sec:ex}.}
A different perspective has been approached recently in two papers by Albano, Cannarsa and Scarinci \cite{AlCaSc2, AlCaSc1}, where they show that if a family of smooth vector fields satisfies the H\"ormander condition, then the set where the local Lipschitz continuity of the minimum time function fails is the union of singular trajectories, and that it is analytic except in a subset on null measure.
{ Finally we mention the recent papers of Bramanti et al. (see \cite{BBMP} and the references therein) on second order 
nonsmooth H\"ormander operators 
where the vector fields  fail to be $C^\infty$ but have the minimal regularity for 
 continuous classical brackets.}

We will proceed in Section 2 with some preliminaries on families of vector fields, in particular the definition of multivalued iterated Lie brackets and deriving the necessary estimates for the corresponding trajectories. In Section 3 we prove the H\"older regularity of the minimum time function for a general nonlinear system in optimal control relaxing the regularity requests on the family of vector fields. In section 4 we turn to H\"older regularity of the solution of the boundary value problem for the eikonal equation and provide some examples where our result is applied to vector fields lacking the classical regularity. Finally in Section 5 we show that our assumptions, under usual 
smoothness of data, are essentially necessary in the case of Lie brackets of length 1.

\section{
Iterated Lie brackets, set-valued extensions, and asymptotic formulas}
This section presents the necessary preliminary definitions and results on Lie brackets for non-smooth vector fields.
\subsection{Classes of 
regularity for Lie brackets} 
In this section we introduce some terminology. 
A vector field $f$ in $\R^n$ is said to be \emph{of class $C^{k, 1}$}  (around a point $x_0\in \R^n$) if it is of class $C^k$ (around $x_0$)
 and its $k$-th order derivatives  are locally Lipschitz continuous (in a neighborhood $V$ of $x_0$). { We use the notation $Df(x)$ for the Jacobian matrix of $f$ at $x$.}
 
 For set-valued vector fields $f: \R^n \ni x \mapsto f(x) \subset \R^n$ 
 we say that $f$ is \emph{of class $C^{-1, 1}$} if it is upper semicontinuous as a set-valued map with compact, convex, nonempty values. 
  
Given  vector fields $f_1, f_2, f_3, f_4, \ldots $ on $\R^N$,  we may compute iterated Lie brackets 
\[
[[f_1, f_2], f_3], \; [[f_1,f_2], [f_3, f_4]], \; [[[f_1, f_2], f_3], f_4], \; [[[f_1,f_2], [f_3, f_4]]  , f_5], \ldots\, ,  
\]
provided that the given vector fields are sufficiently smooth. 
More generally, we may denote any such iterated bracket by $B({ \bf f})$, where ${\bf f} = (f_1, \dots, f_m)$ is a $m$-tuple of vector fields involved in the  definition of $B({ \bf f})$.  $B$ itself may be thought of as 
a {\it (formal) iterated bracket of length $m$} (as a suitable word in a suitable alphabet 
 \cite{FeRa} 
 ), while $B({ \bf f})$ is the result of applying 
  $B$ to ${\bf f}$. 

  We say that {\it ${\bf f}$ is of class $C^B$}, and write ${\bf f} \in C^B$, if
all the components of ${\bf f}$  are  continuously differentiable as many times as it is necessary 
 to compute $B({\bf f})$ so that  
$B({\bf f})$ turns out to be a continuous vector field. E.g., if $B = [[\cdot, \cdot], \cdot]$, then ${\bf f} =(f_1,f_2, f_3) \in C^B$ if and only if $f_1, f_2 \in C^2$ and $f_3 \in C^1$, so that $[[f_1, f_2], f_3]$ is a well-defined continuous vector field; if $B= [[[\cdot, \cdot] , \cdot], \cdot]$, then ${\bf f} = (f_1, f_2, f_3,f_4) \in C^B$ if and only if 
$f_1,f_2 \in C^3$, $f_3 \in C^2$, $f_4\in C^1$, so that $[[[f_1, f_2] , f_3]], f_4]$ is 
 continuous. 

We say that {\it ${\bf f}$ is of class $C^{B-1,1}$}, and write ${\bf f} \in C^{B-1,1}$, if all the components of ${\bf f}$ possess all differentials up to the order that it is necessary to compute 
$B({\bf f})$ minus one, but their highest order differentials are locally Lipschitz continuous; so by virtue of Rademacher's theorem, 
$B({\bf f})(x)$ is well-defined at least for almost every $x \in \R^n$.  E.g., if $B = [[\cdot, \cdot], \cdot]$, then ${\bf f} =(f_1,f_2, f_3) \in C^{B-1,1}$ if and only if
$f_1, f_2 \in C^{1,1}$
 and $f_3 \in C^{0,1}$; if $B= [[[\cdot, \cdot] , \cdot], \cdot]$, then ${\bf f} = (f_1, f_2, f_3,f_4) \in C^{B-1,1}$ if and only if 
$f_1,f_2 \in C^{2, 1}$, $f_3 \in C^{1, 1}$, $f_4\in C^{0,1}$.

\subsection{Multi-flows associated with iterated 
 brackets}
Let us recall that for a (possibly set-valued) vector field $f$ in $\R^n$, 
the flow  $\psi^f$ is the 
 possibly partially defined and possibly set-valued map $\R^n \times \R \ni (x, t) \mapsto \psi^f(x,t) \subset \R^n$ such that for all $(x,t)\in \R^n\times \R$,  $\psi(x,t)$ is the (possibly empty) set of those states $y\in \R^n$ such that there exists an absolutely continuous curve 
$\xi : I_t\to \R^n$ such that $\xi(0) =x$, $\xi(t)=y$, ${\dot\xi}(s) \in f(\xi(s) )$ for a.e. $s\in I_t$, where $I_t = [0\wedge t, 0\vee t ]$. 
The curve $\xi$ is called an \emph{integral curve} of $f$. If $f$ is of class $C^{-1, 1}$, then for any compact $K\subset \R^n$ there exists $T>0$ such that $\psi^f(x, t)$ is not empty for all $(x,t) \in K \times [0, T]$, see \cite{AuCe84}.
Let us call ${\rm Dom}(\psi^f)$, the set of those $(x,t)$ such that $\psi^f(x, t) \neq \emptyset$.  
If $f$ is of class $C^{0,1}$, $\psi^f(x, t)$ is a singleton for each $(x,t)\in {\rm Dom}(\psi^f)$, and we 
view it as a possibly partially defined single-valued map $\R^n\times \R \ni (x,t ) \mapsto \psi^f(x,t) \in \R^n$.

With each formal iterated Lie bracket $B$ of length $m$ and $m$-uple of vector fields ${\bf f}=(f_1, \ldots, f_m)$ of class $C^{B-1, 1}$, we associate 
a  family of multi-flows $\psi_B^{{\bf f}}(t_1, \ldots, t_m)$ for $t_1, \ldots, t_m \in \R$. 
The definition of $\psi_B^{{\bf f}}(t_1, \ldots, t_m)$ is recursive: 
\begin{itemize}
\item if $B$ is a bracket of length $m=1$ so that  ${\bf f}$ consists of a single vector field $f$,
 we set $\psi_B^{{\bf f}}(t) =\psi^f(t)$ for all $t\in\R$, where $\psi^f(t)$ stands for the map $x\mapsto \psi^f(x,t)$;

\item
if $B({\bf f}) = [B_1({\bf f}_{(1)} ), B_2({\bf f}_{(2)} )  ]$, 
where ${\bf f} = (f_1, \dots, f_m)$, ${\bf f}_{(1)} = (f_1, \dots,f_{m_1} )$, ${\bf f}_{(2)} = (f_{m_1+1}, \dots, f_m)$, for $1\le m_1 <m$, we define 
\begin{multline*}
\psi_B^{{\bf f}}(t_1,\dots, t_m) := \psi_{B_2}^{{\bf f}_{2}}(t_{m_1 +1},\dots, t_m)^{-1} \circ \psi_{B_1}^{{\bf f}_{1}}(t_1,\dots, t_{m_1})^{-1} 
 \\ \circ \psi_{B_2}^{{\bf f}_{2}}(t_{m_1 +1},\dots, t_m) \circ \psi_{B_1}^{{\bf f}_{1}}(t_1,\dots, t_{m_1} )
\end{multline*}
for $t_1, \dots, t_m \in \R$.
\end{itemize}
Note that, for $m \ge 2$, the fact that ${\bf f} \in C^{B-1, 1}$ implies that $f_i \in C^{0,1}$ for all $1\le i \le m$, and therefore  
\[
\R^n\times  \underbrace{\R \times \cdots \times \R}_{m-\text{times} } \ni (x, t_1, \ldots, t_m) 
\mapsto \psi_B^{{\bf f}}(t_1, \ldots, t_m)(x) \in \R^n
\]
is a possibly partially defined single-valued map, because all the vector fields are at least locally Lipschitz, and its domain is a nonempty subset of   
$\R^n\times \R^m$. 

For $m=1$, in which case ${\bf f}$ consists just of a vector field $f$ of class $C^{-1,1}$, the (possibly set-valued) map $\psi_B^{{\bf f}}$ is 
 the  flow of $f$. 
  
 \subsection{{ Set-valued iterated 
 brackets} }
Here we give the definition of the set-valued iterated Lie bracket  of length $m$ 
\[
\R^N \ni x \mapsto B({\bf f})_{set}(x) \subset \R^N ,
\]
for a family ${\bf f} =(f_1, \ldots, f_m)$ of $C^{B-1,1}$-regular vector fields. For simplicity we limit 
ourselves to brackets of length $m\le 3$. 
  For iterated brackets of higher length ($m \ge 4$) the reader is referred to \cite{FeRa1}.

 The case of 
 length $m=2$
 is due to   
Rampazzo and Sussmann \cite{RaSu01, RaSu07}:
   for $f_1, f_2 \in C^{0, 1}$ 
 \begin{multline*}
 [f_1, f_2]_{set} (x) := \overline{{\rm co}}\, \big\{ v \; : \; \exists \{x_k\}_{k\in \N} 
 \subset \Diff(f_1) \cap \Diff(f_2)  \\
   \mbox{ such that } x_k \to x \mbox{ as } k \to \infty \; \text{ and }  \; v=\lim_{k\to \infty} [f_1, f_2](x_k)\;  \big\}, 
 \end{multline*}
 where $\Diff(f_i)$ is the set of differentiability points of $f_i$ -- a full measure set by Rademacher's theorem.    
 
It turns out that a mere iteration of this construction to define higher length iterated brackets is  not appropriate
for the validity of the asymptotic formula~{\eqref{asyform}  below}; see \S7 of \cite{RaSu07} for a counterexample.   
 
An appropriate definition for length $m=3$ is the following: 
if $B = [[\cdot, \cdot], \cdot]$ and  ${\bf f} =(f_1,f_2, f_3) \in C^{B-1,1}$, that is, $f_1, f_2 \in C^{1,1}$, $f_3\in C^{0,1}$, one sets
 \begin{multline*}
 [[f_1, f_2],\, f_3]_{set}(x) := \overline{co}\Big\{ v  \; :  \;  
 \exists \{x_k\} \subset \Diff^2(f_1) \cap \Diff^2(f_2)\,, \; \exists \{y_k\} \subset \Diff(f_3)\,, \;\\   
 \mbox{such that } x_k \stackrel{k}{\to } x\,, \; y_k \stackrel{k}{\to }  x  \; \text{ and } \;  
 v=\lim_{k \to \infty} \big( Df_3(y_k) [f_1, f_2](x_k) - D[f_1, f_2](x_k){ f_3}(y_k) \big) \Big\},
\end{multline*}    
    where $\Diff^2(f)$ is the set of points where a vector field $f$ of class $C^{1,1}$ is twice differentiable, a full measure set by Rademacher's theorem.  
    
    $B({\bf f})_{set}(x)$ has 
 nonempty, compact, convex values, it 
  is upper semicontinuous and such that 
$B({\bf f})_{set}(x)$ reduces 
 to the singleton $\{B({\bf f})(x)\}$ at those points $x\in \R^N$ where ${\bf f}$ is of class $C^B$. Hence often  we
 write $B({\bf f})(x)$  instead of $B({\bf f})_{set}(x)$.

\subsection{{ Asymptotic expansions of trajectories} }
\begin{lem}
If  $B$ is a (formal) iterated Lie bracket of length $m$,  ${\bf f} =(f_1, \ldots, f_m)$ is of class $C^B$, $x_*\in \R^N$, then 
\begin{equation}\label{asyform}
\Psi_B^{{\bf f}}(t_1, \dots, t_m)(x) = x + t_1 \cdots t_m  B({\bf f})(x_*) + t_1 \cdots t_m  o(1 )
\end{equation}
as $|(t_1, \ldots,  t_m)| + |x-x_*| \to 0$.
\end{lem}
 This result is classical  for 
  smooth vector fields  (an application of Taylor's formula).  
  Under the minimal 
   regularity assumptions { stated above} it can be found 
   in \cite{FeRa}.
  
{ The regularity assumptions on the vector fields can be further reduced by means of the set-valued brackets $B({\bf f})_{set}(x)$.}
\begin{lem}[{ \cite{RaSu01,RaSu07}  for $m=2$, \cite{FeRa1} for $m\ge 3$.}]\label{as-th}
Given an iterated bracket $B$ of length $m$, ${\bf f }$ of class $C^{B-1,1}(\R^n)$ and $x_*\in \R^n$, then 
\bel{as-for}
{\rm dist}\, \big(\psi_B^{{\bf f}}(t_1, \ldots, t_m)(x) -x, \, t_1 \cdots t_m B({\bf f})(x_*)   \big) = |t_1\cdots t_m| o( 1)
\eeq
as $|(t_1, \ldots, t_m)| +|x-x_*|  \to 0$.

In particular, 
\[
{\rm dist}\, \big(\psi_B^{{\bf f}}(t_1, \ldots, t_m)(x_*) -x_*, \, t_1 \cdots t_m B({\bf f})(x_*)   \big) = o(t_1\cdots t_m)
\quad \text{as $|(t_1,  \ldots, t_m)|   \to 0$.}
\]
\end{lem}
\begin{rem}\label{as-re}
If $m =1$ above, so that  ${\bf f}$ consists of a single vector field $f$ of class $C^{-1,1}$ (the flow $\psi^f$ of $f$ can now be a set-valued map),  then \eqref{as-for} has the following meaning  
\bel{B-exp-m1}
\sup_{y \in \psi^f(x,t) } {\rm dist}\, \big(y-x, \, t f{(x_*)} \big) = |t|\, \gamma(O(|t| + |x-x_*| ) )  
\eeq
as $|t| + |x-x_*| \to 0$, where 
\begin{equation}\label{usc-mod} 
\gamma(\r) = \sup\big\{ {\rm dist}\, \big(w, \, f(x_*)  \big)  \; :\; |x-x_*| \le \r, \;  w\in f(x) \big\}\, .
\end{equation} 
Note that since $f$ is upper semicontinuous at $x_*$, $\g(0+) = \lim_{\r \to 0^+ } \g(\r) =0$, and we call $\gamma$ an \emph{upper semicontinuity modulus}
of $f$ at $x_*$.  
{ The proof of the estimate \eqref{B-exp-m1} is in the Appendix \ref{append}.}
\end{rem}

\section{Sufficient conditions for the H\"older continuity of the minimum time function} 
\label{sufficient}

By a \emph{control system} we mean a family $\mF$ of vector fields on a differential manifold; here for simplicity we limit ourselves to euclidean spaces $\R^n$, for $n\in \N$. See also the following Remark~\ref{old-sys}. 
Let $\mF$ be a control system on $\R^n$. 
By an $\mF$-trajectory we mean any curve obtained as a concatenation of a finite number of integral curves of vector fields in $\mF$.
We say that a control system $\mF$ is symmetric if $-\mF \subset \mF$, where $-\mF = \{-f\; : \; f\in \mF \}$, or more geometrically, any $\mF$-trajectory 
run backward in time is also an $\mF$-trajectory. 

We say that a control system $\mF$ is \emph{(locally) Lipschitz continuous}, or \emph{of class $C^k$}, or \emph{of class $C^{k,1}$}, 
 if any vector field in $\mF$ has such property. 
A system of set-valued vector fields $\mF$ is \emph{of class} $C^{-1,1}$, if each $f\in \mF$ is of class  $C^{-1,1}$.
A \emph{uniformly (locally) 
Lipschitz continuous} control system is a control system $\mF$ such that (for any bounded set $K \subset \R^n$) 
there exists $L \ge 0$ such that $L$ is a Lipschitz constant of $f$ (on $K$)  
for all $f\in \mF$.  
A control system $\mF$ is \emph{uniformly linearly bounded} if there exists $C\ge 0$ such that 
$|f(x)| \le C(|x|+1)$ for all $f\in \mF$, $x\in \R^n$.

Let $x_* \in \R^n$ and $t\ge 0$. The \emph{reachable set} of $\mF$ from $x_*$ at time $t$ is 
\begin{align*}
\mR(x_*, t) &=  \Big\{ y(t)\; : \; y(\cdot) \mbox{ is an } \mF-\mbox{trajectory such that }y(0) =x_* \Big\} \\
&= \Big\{ \psi_{s_1}^{f_1} \circ \cdots \circ \psi_{s_m}^{f_m}(x_*) \; : \; m\in \N, \; i=1, \ldots, m,\;  f_i \in \mF, \; s_i \ge 0, \; \sum_{i=1}^ms_i =t  \Big\} \,
\end{align*}
where $\psi^f_t$, for each $t\in \R$, denotes the map $\R^n \ni x \mapsto \psi^f(x, t) \in \R^n$.  

We say that $\mF$ is \emph{small time locally controllable} (STLC)  from $x_*$ if $x_*$ is an interior point of $\mR(x_*, t)$ for all $t>0$. 
We say that $\mF$ is \emph{(globally) controllable} if for all $x, y \in \R^n$ there exists an $\mF$-trajectory
starting at $x$ and terminating at $y$. 
For symmetric systems we have the following generalization of a classical result of Chow and Rashevski proved in \cite{FeRa2}. 
\begin{defn}[Nonsmooth H\"ormander's condition]\label{Hor}
Let $\mF$ be a control system 
 in $\R^n$ 
  and $x_*\in \R^n$. 
We say that $\mF$ satisfies the {\em nonsmooth H\"ormander's condition}, or the {\em nonsmooth Lie algebra rank condition (LARC)}, at $x_*$, if there exist  formal iterated Lie brackets $B_1, \ldots, B_n$ and tuples of vector fields ${\bf f}_1, \ldots, {\bf f}_n$ of elements of $\mF$ such that
\begin{itemize}
\item ${\bf f}_i$ is of class $C^{B_i-1,1}$ around $x_*$ for all $i=1, \ldots, n$,
\item for all $v_i \in B({\bf f}_i  )(x_*)$, $i=1, \ldots, n$,
\[ {\rm span}\{v_1, \ldots, v_n \}  = \R^n\,. \]
\end{itemize}   
Sometimes the highest length $k\in \N$ of the brackets $B_i$ 
 is relevant, and one says that $\mF$ satisfies 
  {\em H\"ormander's condition of step} $k$ at $x_*$.  

Let $\Omega$ be a subset of $\R^n$.  
One says that $\mF$ satisfies the {\em nonsmooth H\"ormander's condition} (of step $k$), or the nonsmooth LARC,  in $\Omega$  if the property holds 
 at any $x_* \in \Omega$.
\end{defn}
\begin{thm}[A nonsmooth Chow-Rashevski's theorem, \cite{RaSu01}, \cite{FeRa2}]\label{Chow}
Let $\mF$ be a symmetric control system in $\R^n$. 
\begin{enumerate}[(i)]
\item Let $x_*\in \R^n$. If $\mF$ satisfies the (nonsmooth) H\"ormander's condition at $x_*$ of step $k$ for some $k$, 
then
$\mF$ is STLC from $x_*$; moreover, the minimum time function $\R^n \ni x\mapsto T(x,x_*) $
, where 
\bel{mt-pt} 
T(x,x_*) := \inf\{ t\ge 0\;:\; \exists \mF-\mbox{trajectory} \; \xi(\cdot) \; \mbox{such that } \xi(0)=x_*,\; \xi(t)= x   \} 
\eeq
satisfies, for some $C\ge 0$, the estimate 
\[
T(x, x_*) \le C|x-x_*|^{1/k}
\]
in a neighborhood of $x_*$. 
\item If $\Omega$ is an open and connected subset of $\R^n$, and $\mF$ satisfies (the nonsmooth) H\"ormander's condition in $\Omega$ (of step $k$), then 
$\mF$ is globally controllable and the minimum time function is locally H\"older continuous (of exponent $1/k$). 
\end{enumerate}      
\end{thm}

{ Throughout the paper $\mT$ is a closed subset of} $\R^n$ which we shall interpret as a \emph{target} of a control system $\mF$. 
The \emph{minimum time} of $\mF$ to reach $\mT$ is 
$$T(x) := \inf\left\{ t \ge 0\; : \; \mbox{there exists an  $\mF$-trajectory } y(\cdot)  \mbox{ such that }y(0)=x,\; y(t) \in \mT \right\}.$$
The \emph{set of points controllable by $\mF$  to $\mT$} is 
\bel{R-def}
\mR := \left\{ x\in \R^n \; : \; T(x) < \infty \right\}  .
\eeq

\begin{rem}\label{old-sys}
A control system is often given in the form
\bel{sys-old}
\begin{cases}
{\dot y}(t) = f(y(t), \a(t)) \quad t> 0\,, \\
y(0) =x,  \quad x\in \Omega  \,,
\end{cases}
\eeq
where the \emph{control} $\a$ takes values in a given \emph{control set} $A$, a metric space, the state space $\Omega$ is an open subset of $\R^n$, and $f:\Omega \times A\to \R^n$ is continuous and Lipschitz continuous in the state variable, uniformly in the control variable. Clearly a control system in this form can be seen as a control system in the form introduced earlier by considering 
{ 
$\mF^f = \{ f(\cdot, \a)\: :\; \a \in A \}. $ 
The notion of trajectory for 
\eqref{sys-old} 
usually 
admits as \emph{admissible control} any measurable map $\a(\cdot ): [0, +\infty 
] \to A$; 
let us call   \emph{$f$-trajectory} the corresponding \emph{admissible trajectory} of the system~\eqref{sys-old}.
Clearly any $\mF^f$-trajectory is also an $f$-trajectory, being obtained by a piecewise constant admissible control. Although the converse is 
 not true, the former set of trajectories is dense in the latter under rather general conditions.
For a given target $\mT \subset \Omega$, we can 
 define for system \eqref{sys-old} a minimum time function ${\bar T}$  to reach 
 $\mT$ by $f$-trajectories. 
Then ${\bar T}\le T$. 
When $T$ attains continuously the value 0 at $\partial\mT$, it can be shown that in fact ${\bar T} =T$, either by the density property mentioned above or by comparison principles for the Hamilton-Jacobi-Bellman equation satisfied by both $\bar T$ and $T$, see \cite{BaCa97}.
}
\end{rem}

  Given a set $K\subset \R^n$, for $n\in \N$, we denote by $I(K)$ its set of isolated points. 
If $\mT$ is a closed set of $\R^n$--to be interpreted as a target of a control system--throughout this section we use the notation 
\[
d(x) = {\rm dist}\, (x, \mT)
\]
for all $x\in \R^n$.

Our first regularity result for $T$ concerns the case  $\partial\mT$ splits into a $C^1$ manifold and some isolated points. { It is proved by estimating the decrease of the distance function $d$ along admissible trajectories generated by a (possibly non-smooth) Lie bracket via the asymptotic expansion of Lemma \ref{as-th}.}


\begin{thm}
[H\"older continuity of the minimum time - 1]
\label{Hold}
Let $\mF$ be 
symmetric 
and $\mT$ 
such that $I(\mT)$ is locally finite and $\mT\backslash I(\mT)$ is the closure of a nonempty  open set.
 
(i) Let $x_0 \in \pT$.   If  either 

\noindent {(a)}  $\partial\mT$ is of class $C^1$ around $x_0$, and there exist $f_1, \ldots , f_m \in \mF$ 
and a formal iterated Lie bracket $B$ of length $m$ such that 
${\bf f} = (f_1, \dots, f_m)$ is of class $C^{B-1,1}$ in a neighborhood of $x_0$ with
\[
0\notin  B({\bf f})(x_0) \cdot n(x_0) , 
\] 
or (b) $x_0\in I(\pT)$ and $\mF$ satisfies a nonsmooth LARC at $x_0$ of step $m$, 
then 
\[
T(x) \le C\, d(x)^{1/m} 
\]
for some constant $C \ge 0$, in a neighborhood of $x_0$

(ii) Assume that $\mF$ is in addition uniformly locally Lipschitz continuous and uniformly linearly bounded, and $\partial \mT\setminus I(\mT)$ is of class $C^1$. 
If for any $x \in \pT$ condition (i) holds,  possibly with different $f_j$ and $B$, 
 then $\mR$ is open, $T$ is locally H\"older continuous  
 on $\mathcal{R}$, and $\lim_{x\to x_0} T(x) = +\infty$ for all $x_0 \in \partial \mR$.
 
 \noindent If in addition the length of the brackets 
 is at most  $k$ for all $x \in \pT$, 
 then $T$ is locally $(1/k)$-H\"older continuous
 in $\mathcal{R}$.     
\end{thm}
{\bf Proof.} 
The proof of (i) under assumption (b) is just a restatement of Theorem~\ref{Chow}~(i). Indeed, by that 
theorem, for some $C\ge 0$,  $T(x,x_0)\le C|x-x_0|^{1/m}$ for $x$ in a neighborhood  $V$ of $x_0$. Since $x_0\in I(\pT)=I(\mT)$, we can pick $V$ so that in addition $\mT \cap V = \{x_0\}$. On the other hand    
we can find another neighborhood $W\subset V$ of $x_0$ such that any $x\in W$ has a closest point to $\mT$ in  
$\mT \cap V=\{x_0\}$. Therefore, $d(x)= |x-x_0|$ for all $x\in W$, and thus $T(x) \le T(x,x_0)\le C\, d(x)^{1/m}$ for all $x\in W$. 
 
Assume now that $(a)$ holds. The fact that  $\mT$ is of class $C^1$ around $x_0$ 
means by definition that, up to an isometric change of coordinates, $\mT$ is the subgraph of some $C^1$ function in a neighborhood of $x_0$: more precisely, writing $x_0= ( {\bar x}_0, x^n_0) \in \R^{n-1}\times \R$, 
there exists $\varphi: {\bar V}\subset \R^{n-1} \to \R$ of class $C^1$ defined on a closed neighborhood ${\bar V}$ of ${\bar x}_0$, and $\r>0$ such that 
$\mT \cap V  = \{x=({\bar x}, x^n) \in V \;:\:\; x^n \le  \varphi({\bar x}) \}$, where $V={\bar V} \times [-\r, \r]$.  
{ In the change   of coordinates we can also choose the hyperplane of the first $n-1$ coordinates parallel to the hyperplane tangent to $\mT$ at $x_0$, so that $\nabla \varphi({\bar x}_0 ) =0$}. 
 So, if we define the function $w: V\to \R$ by setting $w({\bar x}, x^n)=  x^n - \varphi( {\bar x})$ for $({\bar x}, x^n) \in V$,
$\nabla w(x_0) = {\bf n}(x_0)$, the outer unit normal of $\mT$ at $x_0$, and 
{ a Lipschitz constant of $w$ is $(1+L^2)^{1/2}$ if $L$ is a Lipschitz constant of $\varphi$ on ${\bar V}$. 
Then, since $w(x_0)=0$,
\[
d(x) \le w(x) \le (1+L^2)^{1/2} d(x)
\]
for all $x\in V\cap (\R^n \setminus \mT) $.} 

Let $\xi \in \R^N\cap \mathring{V}$.
By Lemma~\ref{as-for}, for $s>0$, small enough, it is easy to find an $\mF$-trajectory $y(\cdot) : [0, s] \to \R^N$ starting from $\xi$ and satisfying the asymptotic formula
\[
y(s) = \xi + s^m v(s, \xi) + o( s^m ) 
\]
as $ s \to 0$, where $v(s, \xi) \in B({\bf f})(\xi)$ for all $s, \xi$. 
 By the Taylor expansion of $w$ around $\xi$
\[
w(y) = w(\xi) + \nabla w(\xi)\cdot (y-\xi) + o(|y-\xi|),
\]
we find
\[
w(y(s) ) = w(\xi)  + s^m \nabla w (\xi) \cdot v(s, \xi)  + o(s^m )
\]
as $ s \to 0$. 

Since by assumption $0\notin  B({\bf f})(x_0) \cdot \nabla w(x_0)$, it is either 
 $B({\bf f})(x_0) \cdot \nabla w(x_0)\subset ]-\infty, 0[$ or $B({\bf f})(x_0) \cdot \nabla w (x_0)\subset ]0, \infty[$; we may assume that the 
 first holds, changing signs to some of the vector fields if necessary (recall that the system is symmetric). Even more, since  
 $\xi \mapsto B({\bf f})(\xi)$ is an upper semicontinuous set-valued map with compact values and $ \nabla w $ is continuous, for some $\eta >0$, 
 $B({\bf f})(\xi) \cdot \nabla w (\xi)\subset ]-\infty, - 2 \eta[$ for $\xi$ in a neighborhood of $x_0$. 
Therefore  
\[
w(y(s) ) \le w(\xi) -\eta  s^m \le (1+L^2)^{1/2} d(\xi) - \eta  s^m  
\]  
for $\xi $ in a neighborhood of $x_0$, $\xi \notin \mT$, and $s$ small enough, from which it follows
$$d\Big( y \Big(  \big( (1+L^2)^{1/2} d(\xi )/\eta \big)^{1/m} \Big)   \Big)  \le 0.$$ 
It means that the $\mF$-trajectory $y(\cdot)$ has reached the target at a time $s \le \big( (1+L^2)^{1/2} d(\xi )/\eta \big)^{1/m} $, and therefore
$T(\xi) \le   \big( (1+L^2)^{1/2} d(\xi )/\eta \big)^{1/m} $ in a neighborhood of $x_0$. 
   
 (ii) follows from (i)  via the following Lemma~\ref{prop-le} taking  $\o(\r) =  \r^{1/m}$ for $\r\ge0$. 
  \qed 

\medskip
The next lemma is { essentially 
 known in the literature, at least in the case of uniformly globally Lipschitz control system and  bounded $\partial\mT$
 (see, e.g., \cite[Prop. IV.1.6 and Rmk. IV.1.7]{BaCa97} 
 and the  references therein). 
 Here we show how 
a continuity modulus of the minimum time function depends on its 
 modulus at boundary points. }
  
\begin{lem}[Propagation of regularity]\label{prop-le}
Let $\mF$ be a uniformly locally Lipschitz continuous and uniformly { linearly }
 bounded control system (not necessarily symmetric),
$\mR$ the set controllable by $\mF$ to the target $\mT$, and $T: \mR \to [0, \infty[$ the minimum time function. 
Let $\o: [0, \infty[ \to [0, \infty[$ be a 
 modulus 
 such that, for all $x_0 \in \pT$, there is a neighborhood $W$ of $x_0$ and  { $C_0\ge 0$ } satisfying 
  \begin{equation}\label{bnd-reg}
  T(z) \le { C_0}  
    \,\o\big( d( z )  \big)
  \end{equation}
for all $z \in W$.  
Then $\mR$ is open, and   
for any bounded set $V$ with 
$\ol{V}\subseteq\mR$ there exist $C, C_1 \ge 0$ and { $\beta>0$} such that
  \begin{equation}
  \label{holder}
|T(z_1) - T(z_2) | \le C \o (C_1 |z_1 -z_2| ) \quad \forall z_1, z_2 \in V ,  { \;|z_1-z_2|\leq \beta }.
  \end{equation}
Moreover, for all $x_0 \in \partial \mR$, $\lim_{x\to x_0} T(x) =+\infty$.  
\end{lem}
{\bf Proof.}  
Fix $x \in \mR$. 
By the definition of $T(x)$, there exists a 
 sequence $t_k$, 
$T(x)\leq t_k< T(x) + 1/k$, and  $\mF$-trajectories   $y_k(\cdot) : [0, +\infty) 
 \to \R^n$ such
that 
$x=y_k(0),\;x_k = y_k(t_k) \in \partial\mT$.
Assume 
 $y_k(\cdot)$ is 
 the concatenation of integral curves of the 
  vector fields $f_1^k, \ldots, f_{m_k}^k \in \mF$, i.e., 
for suitable times $0=s_0^k<s_1^k <\cdots < s_{m_k}^k  < s_{m_k+1}^k=\infty$,
 the 
maps $\psi_k : \R^n\times \R\ni (z, t) \mapsto \psi_k(z, t) \in \R^n $ defined by 
\begin{equation*}
\psi_k (z
, t) := \psi_{t}^{f_i^k} \circ \psi_{s_{i-1}^k}^{f^k_{i-1}} \circ \cdots \circ \psi_{s_1^k}^{f_1^k}(z
)
\end{equation*}
   if $s_{i-1}^k \le t < s_i^k$, 
$i=1, \ldots, m_k+1$, 
satisfy $y_k(t) =\psi_k(x,t)$ for $t\ge 0$. 
By the assumptions on  $\mF$, 
for any bounded 
$V\subseteq
\R^n$   we have
\begin{gather}
| \psi_k(z, t)  | \le (|z| +C) e^{Ct} \,   \quad \mbox{for all } z \in \R^n,\; t\ge 0\,; \label{Gr-bnd} \\
| \psi_k(z, t) - \psi_k(y, t) | \le  e^{Lt} |z-y|    \quad \mbox{for all } z, y \in V,\; 0\leq t  { \leq   \tau} 
 \,,   \label{Gr-Lip}
\end{gather}
for some $C
 \ge 0$ depending only on $V, \mF$ and  { $ L
 \ge 0$ that may depend also on $\tau$}. 
 
{ Since $t_k\leq T(x)+1$, by estimate \eqref{Gr-bnd} for some $R>0$ $|x_k|\leq R$ for all $k$.
By a covering argument we can find $0<\d \leq 1$ and $C\geq 0$ such that 
\begin{equation}\label{bnd-d-reg}
T(z) \le C \o( d(z) ) \qquad \mbox{for all } z \in B_{\d}(\mT) \cap B_{R+1}(0),  
\end{equation}
 where $B_{\d}(\mT) := \{ z\in \R^n\; : \; d(z) \le \d \}$, $B_{R+1}(0):=\{ z\in \R^n\; : \; |z| \le R+1 \}$.
 
Now 
let $V$ be a bounded neighborhood of $x$, 
take $L$ such that estimate \eqref{Gr-Lip} holds with $\tau= T(x)+1$,
and pick $z\in V$ such that $|z-x|\leq \d e^{-L\tau
}=:\e$.
Then 
$$ d\left(\psi_k(z, t_k), \mT\right)  \le |\psi_k(z, t_k) - \psi_k(x, t_k)  | 
\le 
\d ,$$
and so  $\psi_k(z, t_k) \in B_{\d}(\mT) \cap B_{R+1}(0)$.  
Therefore,  if $T(z) \ge t_k$, we can 
use  the dynamic programming principle
and  \eqref{bnd-d-reg} to estimate
 \begin{align*}
 T(z) -t_k &  \le T\big( \psi_k(z, t_k)   \big)  \le C \o \big( d \big(  \psi_k(z, t_k )    \big)   \big) 
                  \le   C \o  \big( | \psi_k(z, t_k ) -x_k  | \big)\\& =  C \o  \big( | \psi_k(z, t_k ) -\psi_k(x, t_k)   | \big) 
                   \le C \o( e^{L\tau} |z-x| )\,   .
 \end{align*} 
Then $T(z)<+\infty$, and this hold also if $T(z) \leq t_k$, so $\mR$ is open. Moreover, 
  letting $k\to \infty$ we get 
  $T(z) - T(x) \le  C \o( e^{L\tau} |z-x| )$ for all $x, z\in V$, $|z-x|\leq \e$. 
  By exchanging the roles of $x$ and $z$ we obtain the continuity of $T$ at any $x\in\mR$,  and thus $T$ is bounded on any bounded $V$ with
 $\ol{V}\subset \mR$. 
 
Now we take $\tau_V=\max_{\ol{V}}T$ and call $L_V$ the corresponding constant $L$ in \eqref{Gr-Lip}. Then we have  $|T(z) - T(x) |\le  C \o( e^{L_V\tau_V} |z-x| )$ for all $x, z\in V$ such that $|z-x|\leq \beta:=\d e^{-L_V\tau_V}$, which proves \eqref{holder} for $C_1:=e^{L_V\tau_V}$. }

Finally, the proof that $\lim_{x\to x_0} T(x) = \infty$ for all $x_0 \in \partial \mR$ { is the same as}  
 in  \cite{BaCa97}.
  \qed 

\smallskip
Next we extend Theorem~\ref{Hold} 
 to targets 
  with Lipschitz 
  boundaries and  not necessarily $C^1$, or, more generally, 
to targets that satisfy a suitable inner cone condition that we now define. 

For any angle $\theta \in ]0, \pi]$ and unit vector ${\bf n}$, we consider the open unbounded cone with vertex at $0$, opening $\theta$, and axis pointing toward 
   ${\bf n}$: 
 \[
 \Cc_{\theta}( {\bf n} ) = \{h \in \R^n \; : \; h \cdot  {\bf n} > |h| \cos(\theta/2)   \};
 \]   
moreover, for any $\r > 0$, we consider also the 
{ bounded cone} with vertex at $0$, radius $\r$, opening $\theta$ and axis pointing toward ${\bf n}$:
\[
 \Cc_{\r, \theta}( {\bf n} ) = \Cc_{\theta}( {\bf n} ) \cap B_\r(0). 
 \] 

\begin{defn}[Inner cone condition]
\label{cc-def}
Let $\theta \in ]0, \pi]$, and ${\bf n}\in \R^n$ a unit vector.

\noindent For $x\in \pT$ we say that $\mT$ satisfies a
$({\bf n}, \theta)$ {\it inner cone condition} at $x$ if there exists some $\r>0$     
such that $x+ \Cc_{\r, \theta}( {\bf n} ) \subset \mT$.

\noindent For $\Gamma \subset \pT$ we say that $\mT$ satisfies a $({\bf n}, \theta)$ {\it inner cone condition} on 
$\Gamma$ if for all $x\in \Gamma$, $\mT$ satisfies a
$({\bf n}, \theta)$  inner cone condition at $x$. 

\noindent We say that $\mT$ satisfies the {\it inner cone condition}, if for all $x \in \pT$ there exists a neighborhood $\Gamma$ of $x$ in $\pT$, $\theta$, 
and  
${\bf n}$, 
 such that  $\mT$ satisfies a $({\bf n}, \theta)$ inner cone condition on 
$\Gamma$.
\end{defn}

\begin{thm}[H\"older continuity of minimum time -  2]
\label{Hold-bis}
{ Let  $\mF$ be symmetric. 

(i) Let $x_0 \in \pT$ and assume that either 

\noindent { (a)} $\mT$ satisfies a $({\bf n}, \theta)$ inner cone condition in a neighborhood of
$x_0$, 
there are $f_1, \ldots, f_m\in \mF$ 
and
a formal iterated Lie bracket  $B$ of length $m$ such that ${\bf f} = (f_1, \dots, f_m)$ is of class $C^{B-1,1}$ in a neighborhood of $x_0$, and 
\begin{equation}\label{t-br-cond}
  B({\bf f})(x_0) \subset  \Cc_{\theta}({\bf n}) , 
\end{equation} 
or (b) $x_0\in I(\mT)$ and $\mF$ satisfies the (nonsmooth) LARC of step $m$ at $x_0$.
Then, for some constant $C \ge 0$, in a neighborhood of $x_0$
\[
T(x) \le C\, d(x)^{1/m} .
\]

(ii) Assume that $\mF$ is in addition uniformly locally Lipschitz and uniformly linearly bounded.

\noindent  If for any $x \in \pT$ condition  (i) holds, possibly with different  ${\bf n}, \theta, f_j$, and $B$,  
  then $\mR$ is open, $T$ is locally H\"older continuous 
 on $\mathcal{R}$, and $\lim_{x\to x_0} T(x) = \infty$ for all $x_0 \in \partial \mR$.    

\noindent If the length of the brackets  
is at most  $k$ for all $x \in \pT$, 
 then $T$ is locally $(1/k)$-H\"older continuous 
 on $\mathcal{R}$.
 }
\end{thm}
{\bf Proof.}  
The validity of (i) under assumption (b) is an easy 
corollary of Theorem~\ref{Chow}, as in the proof of Theorem~\ref{Hold}. 
So we assume that (a) holds.   
 Let $V$ be a neighborhood of $x_0$ such that $\mT$ satisfies the $({\bf n}, \theta)$ inner cone condition on $\Gamma = V\cap \pT$.
There exists a neighborhood $W$ of
$x_0$ such that each point $\xi \in W$ has a closest point to $\mT$ in $\Gamma$. 
Let $\xi \in W{ \setminus\mT}$ and let ${\bar \xi}$ be a closest point to $\mT$ of $\xi$  in $\Gamma$. Let $\r >0$ be such that    ${\bar \xi} + \Cc_{\r, \theta}( {\bf n} ) \subset \mT$. Let us denote by $e(\cdot)$ the distance function to the cone ${\bar \xi} + \Cc_{\r, \theta}( {\bf n} )$. Clearly $d(\xi) = e(\xi) = |\xi - {\bar \xi} | $. 
Since ${\bar \xi} + \Cc_{\r, \theta}( {\bf n} )$ is a convex set, the function $e(\cdot)$
is differentiable outside the closure of ${\bar \xi} + \Cc_{\r, \theta}( {\bf n} )$, and in particular at $\xi$.

Since $\xi \to  B({\bf f})(\xi)$ is an upper semicontinuas set-valued map with compact values, there exists  $\theta' \in ]0, \theta[ $ such that $B({\bf f})(\xi) \subset  \Cc_{\theta'}({\bf n})$ for any $\xi \in W \setminus \mT$, provided $W$ is taken small enough.
 { Note that $(\xi -\bar \xi) \cdot w \le 0$ for all $w$ in the closure of  $\Cc_{ \theta}( {\bf n})$ because $\bar \xi$ is also the point of the closure of $\bar \xi + \Cc_{\r, \theta}( {\bf n} )$  closest to $\xi$. 
 Thus $\xi -\bar \xi$ forms an angle $\ge \pi/2$ with any vector in the closure of $\Cc_{ \theta}( {\bf n})$. For every $v\in B({\bf f})(\xi)$ every vector $w$ that forms an angle $(\theta -\theta')/2$ with $v$  belongs to the closure of $\Cc_{ \theta}( {\bf n})$, therefore the angle between the vector $\nabla e(\xi) = (\xi -{\bar \xi})/|\xi -{\bar \xi}|$ and any vector  $v\in B({\bf f})(\xi)$ is greater than $(\theta - \theta' + \pi)/2$, for $\xi \in W\setminus \mT$. Thus $\nabla e(\xi) \cdot v  \le -\sin((\theta - \theta')/2) |v| $ for all $v \in  B({\bf f})(\xi)$, $\xi \in W\setminus \mT$.}

	By Lemma~\ref{as-for} it is easy to find an $\mF$-trajectory $y(\cdot) : [0, s] \to \R^N$ starting from $\xi$ satisfying the  asymptotic formula
\begin{equation}\label{traj-est}
y(s) = \xi + s^m v(s, \xi) + o( s^m )  
\end{equation}
as $ s \to 0$, where $v(s, \xi) \in B({\bf f})(\xi)$ for all $s, \xi$. 
Then  we can estimate as follows 
\begin{multline*}
d( y(s)) - d( \xi ) \le e(y(s) ) - e(\xi) = s^m \nabla e(\xi) \cdot v(s, \xi) + o(s^m) \\
 \le - s^m \sin((\theta - \theta')/2) |v(s, \xi) | +  o(s^m), 
\end{multline*}  
and thus for $s$ small enough
 \[
d( y(s)) \le  d( \xi ) - \eta s^m, 
 \]
where $ \eta =  (1/2) \e \sin((\theta - \theta')/2)$ with $\e>0$ such that $|v| \ge \e$ for all $v \in B({\bf f})(\xi )$, $\xi \in W$. 
Thus, provided we chose $W$ small enough, we have $d\Big( y \Big(  \big( d(\xi )/\eta \big)^{1/m} \Big)   \Big)  \le 0$. 
It means that the trajectory $y(\cdot)$ has reached the target at a time $s \le \big( d(\xi )/\eta \big)^{1/m} $, and therefore
$T(\xi) \le   \big( d(\xi )/\eta \big)^{1/m} $ in a neighborhood of $x_0$. 

Clearly (ii) follows from (i) via Lemma~\ref{prop-le}.
\qed 

{
\begin{rem}[An alternative via the 
 superdifferential of the distance]
	The previous theorem continues to hold if instead of condition (a)  
	we require the following: there exist $R_{x_0}>0$ and a neighborhood $V$ 
	of $x_0$ such that for all $\xi \in V\setminus \mT$ there exist { $\nu(\xi) \in D^+
	 d(\xi)$} (the usual Fr\'echet superdifferential of $d$ at $\xi$), a formal iterated Lie bracket $B$ of length $m$, and vector fields ${\bf f} = (f_1, \dots, f_m)$ in $\mF$ of class $C^{B-1,1}$ on $V$ 
	so that 
	\[
	\nu(\xi) \cdot B({\bf f}) (\xi ) < -R_{x_0}\,, 
	\]
	{ i.e., $\nu(\xi) \cdot v< -R_{x_0}$ for every $v \in B({\bf f}) (\xi )$.}
	The proof is similar: { $\nu(\xi) \in D^+ d(\xi)$ means} 
	\[
	d(y) - d(\xi) \le \nu(\xi) \cdot(y-\xi) +o(y-\xi) \quad\text{as $y\to \xi$} .
	\]
	Combining this  inequality  with the estimate  \eqref{traj-est} for an  $\mF$-trajectory $y(\cdot )\colon[0,s]\to \R^N$ as in the previous proof,
 we obtain 
\[
	d( y(s)) - d( \xi ) \le -s^m \nu(\xi) R_0 +o(s^m)   \quad\text{ as $s\to 0$}
	\]
	 for all $\xi \in W\setminus \mT$, where $W$ is a small enough neighborhood of $x_0$. 
	The rest then follows exactly as in the previous proof. 
	
	{ Note that the assumption of non-empty $D^+ d(\xi)$ at points $\xi\in\partial\mT$ implies more regularity of $\partial\mT$ than in Theorem  \ref{Hold-bis}.}
	\end{rem}
}

\section{Degenerate and nonsmooth eikonal equations} 
\label{eikonal}
\subsection{Problem statement and setting}
In this section we establish results on the solvability, continuity, and especially H\"older regularity, of viscosity solutions to the Dirichlet boundary value problems associated with a large class of degenerate eikonal equations 
with quite nonsmooth  coefficients and nonsmooth domains.  


The Dirichlet problem associated with a typical eikonal equation that we study is
\bel{eik}
\left\{ \begin{aligned}
& \sum_{i=1}^m| f_i u|^2 + 2\sum_{i=1}^m b_i(x)  f_iu = h^2(x)\, ,  \quad \mbox{ in } \Omega\\  
&u =g \qquad  \mbox{ on } \partial \Omega, 
\end{aligned} \right.
\eeq  
where $\Omega$ is some open subset of a differentiable manifold $M$, $f, \ldots, f_m$ are vector fields on $M$, 
 $b_1, \ldots, b_m$, $h$ are functions defined on $M$, and $g$ is a function defined on $\partial \Omega$. 
Clearly, we are using the usual identification of vector fields with first-order partial differential operators: more precisely, if in a coordinate chart 
$f_i(x) = \left(f_i^1(x), \ldots, f_i^n(x) \right)$, where $n= {\rm dim M}$, then 
$f_i u(x) = \sum_{j=1}^n f_i^j(x) \partial_{x_j}u(x)$.

In the literature the eikonal equation is often written  in the following nonintrinsic (coordinate-dependent) form 
\[
|\s(x)^t \nabla u(x)|^2 +2 b(x)\cdot \s(x)^t \nabla u(x) = h^2(x)\,  ,
\]
where $\s(x)$ is the matrix with columns $f_1, \ldots, f_m$ (i.e., their coordinate representations), 
$b(x) = (b_1(x), \ldots, b_m(x))$. 
For simplicity we will work in $M = \R^n$, 
although 
the results stated here could be extended to more general differential manifolds.

Simple computations allow  to rewrite the PDE in \eqref{eik} in the following Hamilton-Jacobi-Bellman form
 \bel{eik-bis}
\left\{
 \begin{aligned}
&\max_{\a \in B'_1(0)  } \left\{  \sum_{i=1}^m \a_i f_i u (x)  - \ell(x, \a) \right\} = 0 & &\mbox{ in } \Omega\subseteq\R^n
\\  
&u =g &  &\mbox{ on } \partial \Omega, 
\end{aligned}
\right.
\eeq        
where 
\bel{ell-def}
\ell(x, \a) = { \left( h(x)^2 + |b(x) |^2  \right)^{1/2} - \sum_{i=1}^m\a_i b_i(x) }
\eeq
for all $x\in  \R^n$, $\a=( \a_1, \ldots, \a_m) \in B'_1(0) $,
and  $B'_1(0)$ is the closed unit ball of $\R^m$.
Then we can give a control interpretation of the problem via the 
 symmetric control system
\bel{cs-eik}
\left\{ \begin{aligned}
&{\dot y} = \sum_{i=1}^m \a_i f_i(y)   \\ 
&y(0) = x,  
\end{aligned} \right.
\eeq
where the set of admissible controls denoted by $\mA$ consists of all the measurable maps $\a(\cdot)  = (\a_1(\cdot) , \ldots, \a_m(\cdot) ) : [0, \infty[ \to 
B'_1(0)$. 
{
This is in the form \eqref{sys-old} described in Remark \ref{old-sys}, and corresponds to the family of vector fields 
\[
  \mF_0 = \left\{ \sum_{i=1}^m \a_i f_i \; : \;  \a = (\a_1, \ldots, \a_m )\in B'_1(0) \right\}  
  \]
For each $x \in \R$ and $ \a(\cdot) \in \mA$ let $t \mapsto y(t; x, \a(\cdot))$ be the solution of \eqref{cs-eik}, 
and denote $\t_x(\a(\cdot) )$  the first time it hits the target $\mT = \R^n \setminus \Omega$. 
Then the minimum time function can be written as $T(x) = \inf\{\t_x(\a(\cdot) ) \,:\, \a(\cdot) \in \mA\}$, and the set $\mathcal{R}$ is defined by \eqref{R-def} in the previous section. 
   A candidate solution of \eqref{eik-bis} is } the value function of the following optimal control problem 
\bel{op-pr}
v(x) := \inf_{\a(\cdot) \in \mathcal{A} } \left\{ \int_0^{\t_x(\a(\cdot) )} \ell \left( y(t; x, \a(\cdot)) , \a(t)  \right) dt + g\left( y(\t_x; x, \a(\cdot))  \right) \right\}.  
\eeq
{ However, this is not true if $\Omega$ is not a subset of the reachable set $\mathcal{R}$, and in that case \eqref{eik-bis} does not have a solution. More precisely, it was proved in \cite{BaSo91} (see also \cite{BaCa97}) that there is at most one open set  $\mathcal{O}\subseteq \Omega$ such that there is a continuous solution $u$ of
\bel{eik-fb}
\left\{ \begin{aligned}
&\sum_{i=1}^m| f_i u|^2 + 2\sum_{i=1}^m b_i(x)  f_i u = h^2(x) \quad \mbox{ in } \mathcal{O}  \cap \Omega ,\\
&u =g \qquad  \mbox{ on } \partial \Omega \\
& u(x) \to \infty  \qquad \mbox{ as } x\to \partial \mathcal{O} . 
\end{aligned} \right.
\eeq  
Moreover, if $v$ is continuous at all points of $ \partial \Omega$, then  the pair $(\mR, v)$ is the unique 
 solution of \eqref{eik-fb}. So there is a solution of \eqref{eik} and  \eqref{eik-bis} if in addition $\Omega\subseteq\mathcal{R}$. The continuity of $v$ depends on the controllability of \eqref{cs-eik} near $\partial\Omega$ and on a compatibility condition on the boundary data $g$, that we now recall.
 }
For all $x \in \overline{\Omega}$, $z \in \partial \Omega$ we define 
\begin{multline}\label{comp-L}
L(x,z) := \inf\Big\{ \int_0^t l (y(s), \a(s) ) ds \; : \; t \ge 0, \; y(\cdot) \in AC([0, t]),\;  y(0) =x,\; y(t) = z,  \\ 
y(s) \in \Omega \; \forall s\in ]0,t[, \; a(\cdot) \in \mA,\; {\dot y}(s) = \s( y(s)) \a(s) \mbox{ for a.e. } s \in [0,t]   \Big\}.  
\end{multline}
We will assume that $g$ satisfies the {\em compatibility condition} 
\bel{comp-con}
g(x) - g(z) \le L(x, z) \quad\forall x,z \in \partial\Omega. 
\eeq
Sufficient conditions for this compatibility condition to hold are
\bel{comp-suff-1}
g(x) - g(z) \le 
\frac1{L_o} \ln\left(1 + \frac{L_o|x-z|}{ C}  \right) \, ,
\eeq
or 
\bel{comp-suff-2} 
g(x) - g(z) \le C^{-1} |x-z|  
\eeq
for all $x,z\in \Omega$,
where $C$ and $L_o$ are, respectively, an upper bound for $|f_i|$ and the Lipschitz constants of $f_i$ for all $i=1,\dots,m$,  in $\Omega$, see \cite[Proposition IV.3.7]{BaCa97}.

\subsection{Regularity of the solution}

\begin{defn}[H\"ormander's condition at the boundary]
\label{def:Hbound}
Let $x_0 \in \partial {\Omega}$. We say that $\mF = \{f_1, \dots, f_m\}$ (or the matrix-valued function $a= \s\s^t$, or the Dirichlet problem \eqref{eik}) is {\em  H-noncharacteristic of degree 
$k\in\N$ at $x_0$}  if 
 either one of the following holds:
 \noindent
(i) $x_0\in I(\partial \Omega)$ and $\mF$ satisfies (the possibly nonsmooth) LARC 
 of step $k$;
\noindent
(ii) there exist  a unit vector ${\bf n} \in \R^n$ and an angle $\theta \in ]0, \pi]$ such that 
$\mT= \R^n \setminus \Omega$ satisfies the $({\bf n}, \theta)$ inner cone condition in a neighborhood of $x_0$ in $\pT= \partial \Omega$,
 there exists 
 a formal iterated Lie bracket $B$ of length $k$ and $i_1, \ldots, i_k \in \{ 1, \ldots, m\} $ such 
that ${\bf f} = (f_{i_1}, \ldots, f_{i_k})$ is of class $C^{B-1,1}$ in a neighborhood of $x_0$, and 
\bel{entrante}
B( {\bf f} )(x_0) \subset \mC_\theta({\bf n} ).   
\eeq
{ We say that $\mF$ 
 is {\em H-noncharacteristic} at the boundary $\partial \Omega$ 
 if, for every $x_0 \in \partial \Omega$, $\mF$ is  H-noncharacteristic of degree $k$ at $x_0$ for some $k $, and that it 
 is {\em H-noncharacteristic of degree $k
 $} 
 if, for every $x_0 \in \partial \Omega$, $\mF$ is  H-noncharacteristic of degree $k'$ at $x_0$ for some $k' \le k$, and at some point $k'=k$.}
\end{defn}
\begin{rem}
Observe that the condition \eqref{entrante} of Definition (\ref{def:Hbound}) is equivalent to 
\[
 0 \notin B( {\bf f})(x_0) \cdot n(x_0)    
\]
if $ \Omega $ is of class $C^1$ in a neighborhood of $x_0$ and $n(x_0) $ is its outward normal at $x_0$, and to $B( {\bf f})(x_0) \cdot n(x_0)\ne 0$ if the bracket is single-valued.
\end{rem}

\begin{thm}\label{eik-thm}
Let $\Omega\subset \R^n$ 
   $\mT= \R^n\setminus \Omega$ 
satisfies the inner cone condition in $\pT\setminus I(\pT)$, 
 the vector fields in $\mF= \{f_1, \ldots, f_m \}$  be Lipschitz continuous, $b$, $h$, $g$ continuous, { $h(x) \ne 
 0$} for all $x\in \Omega$, and $g$ satisfy the compatibility condition \eqref{comp-con}. 
 Assume also that $\mF$ is H-noncharacteristic at the boundary $\partial \Omega$.  
  Then the following facts hold.  
  
  (i) $\mR$, defined by \eqref{R-def},  is open and contains $\mT = \R^n\setminus \Omega$, 
  $\mR =\{x\in \R^n\; :\; v(x) <\infty\}$,  $v$ is  bounded below and continuous on 
  $\mR$: more precisely, on any bounded set $V\subset \mR$, $v$ has a
  continuity modulus $\omega_{v, V}$ 
   of the form
  \[
  \o_{v, V}(\r) = C \o_\ell( C \r) + C \r^{1/k}  +  \o_g ( C \r^{1/k} +C \r )+ C \r \,,\quad \r \ge 0 ,
  \]
 where $\omega_g$, $\omega_\ell$ 
   are continuity moduli of $g, \ell$ { restricted to a bounded set $K$ 
    depending on  $V$ and the data}.  In particular, if $g, \ell$ are locally H\"older continuous, so is $v$.     
Moreover, 
{ the pair $(\mR, v)$ is the unique 
 solution of \eqref{eik-fb}.}

(ii) If $\mF$ is H-noncharacteristic of degree $k$ 
 at the boundary $\partial \Omega$, $b$, $h$ are locally $(1/k)$-H\"older continuous and $g$ is Lipschitz continuous, then $v$ is 
locally  $(1/k)$-H\"older continuous on $\mR$. 

(iii) If in addition $\mF$ satisfies  the nonsmooth H\"ormander's condition in the interior of $\Omega$, then $\mR = \R^n$ and $v \in C(\R^n)$ 
is the unique  (continuous) viscosity solution of \eqref{eik} bounded from below.   
\end{thm}
\begin{rem}\label{opt-expon}
 We point out that the H\"older continuity exponent $1/k$ in Theorem~\ref{eik-thm}~(ii),~(iv)  is optimal. 
To see this consider the problem studied in \cite{Al12}, see in particular Remark~1.1 (i) of that paper.   
\end{rem}
{\bf Proof of Theorem~\ref{eik-thm}.} The openness of $\mR$ is stated by Theorem~\ref{Hold-bis}~(ii). 
Let $x \in \R^n$.  
{ Note that  $h^2>0$ implies $\ell>0$, and so} $M_\ell T(x) + M_g \ge v(x) \ge m_\ell T(x) - M_g$, where $M_\ell\ge m_\ell>0$ are, respectively, an upper and a lower bound of $\ell$ 
 on the set $\{ y(t; x, \a(\cdot)) \;:\; 0\le t \le T(x), \; \a(\cdot) \in \mA \}$ which is bounded, while $M_g$ is an upper bound of $|g|$. 
Therefore $v(x) \to \infty $ as $x \to \partial \mR$ by Lemma~\ref{prop-le}, and 
$\mR = \{ x\in\R^n \: : \; v(x) <\infty \}$. 

 It is easy to verify that 
\[
v(x) = \inf_{z\in \partial \mT} \big( L(x, z) + g(z) \big)  \quad \forall x \in \ol{\Omega}, 
\]  
where $L$ is defined by \eqref{comp-L}. Note that $L$ satisfies 
 the triangle inequality $L(z, z') \le L(z, x) + L(x, z')$ 
  and symmetry $L(x,z') = L(z',x)$, because 
  the  control system \eqref{cs-eik} is symmetric.  
Therefore, using 
also the compatibility condition \eqref{comp-con}, we have
{
\begin{align*}
g(z) -v(x) &= \sup_{z'\in \partial \mT}\big(g(z) - g(z') -L(x, z')  \big) 
\le L(z, x) \qquad \forall z \in \partial \Omega, x \in \ol{\Omega}. 
\end{align*}
Let $x_0 \in \pT$ and  $W$ a bounded neighborhood of $x_0$ 
 such that 
$T(x) \le C d(x)^{1/k} $ in $W$, 
by Theorem~\ref{Hold-bis}~(i). }
 Let $x\in W$ and $\e>0$. There exists a control $\a(\cdot)$ such that $\t_x(\a(\cdot) ) < T(x) + \min\{\e, 1\}$. 
 The set of points reachable from $W$  up to time $\sup_{x\in W}T(x) +1$, say $K$, is bounded.
 Let $M_\ell$, $M_f$ be bounds of { $\ell (x, \a)$ and $|\sum_{i=1}^m \a_i f_i(x)|$ on $K\times B'_1(0)$,}
 and let $\o_g$ denote the continuity modulus of $g$ on $\pT \cap W$. 
Then,  on one hand
\[
g(x_0) -v(x) \le L(x_0, x) \le M_\ell (T(x) +\e) \le M_\ell (C d(x)^{1/k} +\e)\leq M_\ell(C|x-x_0|^{1/k}+\varepsilon), 
\]
and on the other 
\begin{align*}
v(x) - g(x_0) &\le \int_0^{\t_x(\a(\cdot) )} \ell \big( y(t; x, \a(\cdot)) , \a(t)  \big) dt + g\big( y(\t_x; x, \a(\cdot))  \big) - g(x_0) \\
&\le  M_\ell( T(x) +  \e) + \o_g( |y(\t_x; x, \a(\cdot)) -x  | + |x -x_0| ) \\
&\le  M_\ell (C d(x)^{1/k} +  \e) +  \o_g( M_f (T(x) +\e) + |x -x_0| ) \\
&\le M_\ell ( C |x-x_0|^{1/k} +\e)  +  \o_g ( M_f (C |x-x_0|^{1/k} + \e) + |x -x_0|    ).
\end{align*}
Since $\e>0$ is arbitrary, we have shown that
$|v(x) -g(x_0)| \le \o_{v, W} (|x-x_0|) $, where 
$\o_{v, W}(\r) = M_\ell ( C \r^{1/k})  +  \o_g ( M_f (C \r^{1/k}) + \r )$ for all $\r\ge 0$.  

 Now let $V \subset \mR$ be bounded. { We are going to show 
 that  the continuity of $v$ at boundary points propagates in 
 $V$, with an estimate of 
 the modulus of continuity.} 
 %
 %
 Let $z_1, z_2\in V$ with $|z_1 -z_2| \le {\bar \d}$, where ${\bar \d }>0$ is chosen below. Let also $\e \in ]0, 1]$. There exists a control $\a(\cdot) \in \mA$ such that
\[
v(z_1) +\e \ge \int_0^{\t_{z_1}(\a(\cdot) )} \ell \big( y(t; z_1, \a(\cdot)) , \a(t)  \big) dt + g\big( y(\t_{z_1}; z_1, \a(\cdot))  \big).
\]
Since $v$ is locally bounded 
 and so are $\ell$ 
  and $g$, it follows that 
$\t_{z_1}(\a(\cdot) ) \le \t_0$ for some $\t_0\ge 0$ that depends on  $\Omega, \mF, \ell, g, V$ but not on the particular $z_1\in V$ or $\e \in ]0,1[$. 
 Let $K$ denote the set of points $x$ reachable from $V$ up to time $\t_0$, which is bounded.
It is possible to find $\d>0$, $C\ge 0$, $k \in \N$
 (which depend on $K, \ell, g, \Omega$, and hence on $V, \mF, \ell, g, \Omega$) such 
that $|v(x) -g(x_0)| \le {\bar \o}_{v, V} (|x-z_0|) $ for all $x \in B_\d(\pT) \cap K$, $x_0 \in \pT \cap K$, 
 where 
${\bar \o}_{v, V}(\r) = C \r^{1/k} +  \o_g ( C \r^{1/k} + \r ) +C\r$ for all $\r\ge 0$, $B_\d(\pT) = \{x \in \R^n \; : \; {\rm dist}(x, \pT) \le \d \}$, and
$\o_g$ is the continuity modulus of $g$ on $\pT \cap K$. Since the maps $ z\mapsto y(t; z, \a(\cdot) )$ are (locally) Lipschitz, uniformly in $t \in [0, \t_0]$, $\a(\cdot) \in \mA$, we choose ${\bar \d} >0$ such that
$|z_1 -z_2| \le {\bar \d}$ implies $|y(t; z_1, \a(\cdot) ) -y(t; z_2, \a(\cdot) ) | \le \d$.

Suppose first that $t_2:=\t_{z_2}(\a(\cdot) ) \le \t_{z_1}(\a(\cdot) ):= t_1 $. Then we have 
\begin{align}
v(z_2) -v(z_1) &\le \int_0^{t_2} \ell \big( y(t; z_2, \a(\cdot)) , \a(t)  \big) dt + g\big( y(t_2  ; z_2, \a(\cdot))  \big) \nonumber
 \\ &\qquad\qquad\qquad\quad-  \int_0^{t_1} \ell \big( y(t; z_1, \a(\cdot)) , \a(t)  \big) dt -  g\big( y(t_1; z_1, \a(\cdot))  \big) +\e\nonumber
\\ &= \int_0^{t_2} \Big( \ell \big( y(t; z_2, \a(\cdot)) , \a(t)  \big)  - \ell \big( y(t; z_1, \a(\cdot)) , \a(t)  \big) \Big) dt \nonumber \\ 
  & \; - \int_{0}^{t_1 -t_2} \ell \big( y\big(t ; y(t_2; z_1, \a(\cdot) ), \a(\cdot +t_2 )  \big) , \a(t+t_2)  \big) dt \nonumber\\
  & \quad - g\big( y\big(t_1  - t_2 ; y(t_2; z_1, \a(\cdot) ) , \a(\cdot +t_2) \big) \big)   
  +g\big( y(t_2  ; z_2, \a(\cdot) )  \big) +\e \nonumber\\
&\le \t_0 \o_\ell ( L_y |z_1- z_2|) +  g\big( y(t_2  ; z_2, \a(\cdot) )  \big) - v\big( y(t_2; z_1, \a(\cdot) )  \big) +\e\nonumber \\
&\le \t_0 \o_\ell ( L_y |z_1- z_2|) + { \bar{\o}}_{v, V} ( |y(t_2  ; z_2, \a(\cdot) )  -  y(t_2; z_1, \a(\cdot) ) | ) + \nonumber\e \\
&\le \t_0 \o_\ell ( L_y |z_1- z_2|) + \bar{\o}_{v, V} ( L_y |z_1 -z_2|)+\e, \nonumber
 \end{align}
where $L_y\ge 0$ is a common Lipschitz constant on $K$ of the maps $z \mapsto y(t; z, \a(\cdot))$, for $t\in [0, \t_0]$, $\a(\cdot) \in \mA$, and 
$\o_\ell$ is a common continuity modulus of the maps $x \mapsto \ell(x, \a)$, for $\a \in A$, on $K$;
notice that above we have also used the fact that  $y(t_2  ; z_1, \a(\cdot) )  \in B_\d(\pT) \cap K$: this is true because 
$y(t_2  ; z_2, \a(\cdot) ) \in \pT \cap K $ and $|z_1-z_2| \le {\bar \d}$ implies $|y(t_2  ; z_1, \a(\cdot) ) -y(t_2  ; z_2, \a(\cdot) ) | \le \d$.  

If, instead, $t_2 \ge t_1$ we use 
the dynamic programming principle to obtain
 \begin{align*}
 v(z_2) -v(z_1) &\le \int_0^{t_1} \ell\big(y(t; z_2, \a(\cdot) )   \big)dt + v\big(  y(t_1; z_2. \a(\cdot))  \big) - v(z_1)  \\
 &\le \int_0^{t_1} \ell\big(y(t; z_2, \a(\cdot) )   \big)dt  -\int_0^{t_1} \ell \big( y(t; z_1, \a(\cdot)) , \a(t)  \big) dt  \\ &\quad+ v\big(  y(t_1; z_2. \a(\cdot))  \big)  -  g\big( y(t_1; z_1, \a(\cdot))  \big) +\e ,
 \end{align*}
which is estimated as above. The roles of $z_1, z_2$ can be exchanged, and letting $\e \to 0$, we have
 proved $|v(z_2) -v(z_1)| \le \t_0 \o_\ell ( L_y |z_1- z_2|) + {\bar \o}_{v, V} ( L_y |z_1 -z_2|) $ for all $z_1, z_2 \in V$ with $|z_1- z_2|\le {\bar \d} $.  
 
From this and the boundedness of  $V$ it  follows that a  
  continuity modulus $\o_{v, V }$ of $v$ on $V$ 
  has the form  $\o_{v, V}(\r) = C \o_\ell( C r) + C \r^{1/k}  +  \o_g ( C \r^{1/k} +C \r )+ C \r$, for $\r \ge 0$,  for some 
  $0\le C< \infty $ and $k\in \N$ that depend on $V, f, \ell, g, \Omega$; under assumptions of (ii), the value of $k$  coincides with that of (ii). Therefore we have proved all the claims about the regularity of $v$.   
    
{ Once we know that the value function is continuous, it is standard in viscosity solutions theory that it satisfies \eqref{eik-fb} \cite{BaCa97}. 
The uniqueness of $(\mR, v)$ is 
proved in \cite[Theorem~3.1]{BaSo91}. 
 Finally,  the fact that $\mR = \R^n$ in (iii) }
  follows from the (nonsmooth) Chow-Rashevski's Theorem~\ref{Chow}. 
\qed
  
\subsection{Examples}
\label{sec:ex}

\begin{example}[Nonholonomic integrator, or Brockett's vector fields, or generators of the Heisenberg group] 
In $\R^3$ consider the control system $\mF= \{f_1,f_2\}$ with 
\[
f_1(p) = \left(\begin{array}{c} 1\\ 0\\ -y \end{array} \right) \equiv \partial_x - y \partial_z, \quad   
f_2(p) = \left(\begin{array}{c} 0\\ 1\\ x \end{array} \right) \equiv \partial_y + x \partial_z \,,
\] 
for all $p=(x,y,z)^t\in\R^3$. 
One easily checks that $[f_1, f_2] =2\partial_z$. 
Thus $\mF$ satisfies the LARC 
 of step $2$ at any point of $\R^3$.
The eikonal Dirichlet problem is
\begin{equation}\label{DP-ex-1}
\begin{cases} |\partial_x u -y \partial_z u|^2 + |\partial_y u + x \partial_z u|^2 =1 & \mbox{in $\Omega$} ,\\
                     u =g & \mbox{on $\partial \Omega$ } .  \end{cases}                  
\end{equation} 
For $g \equiv 0$, Theorem~\ref{eik-thm} 
gives a unique locally $(1/2)$-H\"older continuous viscosity solutions on the closure of any open domain $\Omega \subset \R^3$ 
whose complement $\mT = \R^3 \setminus \Omega$ is such that 
$\mT \setminus I(\mT)$ satisfies an inner cone condition, see Definition~\ref{cc-def}~(iii).   
Furthermore, this solution coincides with the minimum time function $T$ to reach the target $\mT$ by trajectories of 
$\mF_0= \{\pm f_1, \pm f_2 \}$.  

For general continuous boundary data $g$, 
 problem~\eqref{DP-ex-1} still admits  a unique continuous viscosity solution on 
$\ol{\Omega}$ 
 provided that $g$
satisfies the compatibility condition \eqref{comp-con}. For instance, \eqref{comp-con} holds if
 $g(p) -g(q) \le C^{-1} |p - q|$ for all $p, q\in \partial\Omega$, with $C = 1 + \max\{ |x| \vee |y| :   (x, y, x) \in \Omega\}$. If $g$ is in addition locally H\"older continuous, then
the solution is also locally H\"older continuous on $\ol{\Omega}$. If $g$ is locally Lipschitz, then the solution is locally $(1/2)$-H\"older continuous on $\ol{\Omega}$.  

{ To arrive at a locally Lipschitz solution of \eqref{DP-ex-1}, 
 we must 
  assume not only 
    $g$  locally Lipschitz and satisfying 
    \eqref{comp-con}, but also that all $p\in 
    \partial \Omega$ are truly noncharacteristic. For $\Omega$ of class $C^1$ this means that  $ {\bf n} (p) \cdot f_i(p) \neq 0$ for some $i$. In our non-smooth context it means that
$\mT$ satisfies a $({\bf n}(p), \theta(p))$ inner cone condition 
on some relative 
neighborhood of $p$ in the boundary $\partial\Omega$, and  $\mC_{\theta(p)}({\bf n}(p) ) \cap span\{f_1(p), f_2(p)\} \neq \{0\}$ .}
 \end{example} 
 \begin{example}[Nonsmooth Brockett type vector fields]
Consider the vector fields $\mF= \{f_1, f_2\}$ on $\R^3$ defined by setting, 
 for  $p= (x,y,z)^t\in \R^3$, 
 \[
f_1(p) = \left(\begin{array}{c} 1\\ 0\\ \a(y) \end{array} \right) \equiv \partial_x + \a(y)  \partial_z, \quad   
f_2(p) = \left(\begin{array}{c} 0\\ 1\\ \b(x) \end{array} \right) \equiv \partial_y + \b(x)  \partial_z \,,
\]
where $\a, \beta : \R\to \R$ are Lipschitz continuous functions. The Lie bracket  $[f_1, f_2]$ can be computed classically 
at the points 
where $\a, \b$ are both differentiable, and 
at those points 
 $[f_1, f_2](p) = (\beta'(x)-\a'(y))\partial_z$. We can 
 compute $[f_1, f_2]$ in the set-valued sense of this paper 
 in terms of Clarke's generalized derivatives $D_C$ of $\a$ and $\beta$ and get
\[
[f_1, f_2]_{set}(p) =  \big(D_C\beta(x) - D_C \a(y) \big)\partial_z\,  
\] 
for all $p=(x,y,z) \in \R^3$. 
For instance, if $\a, \b$ are continuous and piecewise  $C^1$, then for all $p=(x,y,z)^t \in \R^3$    
%
 %
 \[
[f_1, f_2]_{set}(p) = [m(x,y), M(x,y)] \partial_z =\{\lambda \partial_z \;:\; \lambda \in [m(x,y), M(x,y)] \},
\]
where 
\begin{equation*}
\begin{gathered}
m(x, y) = \min\left\{ b-a\; : \; b \in \{ \beta'(x-), \beta'(x+)\}, \;   a\in \{\a'(y-), \a'(y+)\} \right\}, \\  
M(x,y)= \max\{ b-a\; : \; b \in \{ \beta'(x-), \beta'(x+)\}, \;   a\in \{\a'(y-), \a'(y+)\} \} \,. 
\end{gathered}
\end{equation*}
If 
\begin{equation}\label{hor-c-es}
0 \notin D_C\beta(x) - D_C \a(y)  
\end{equation}
or, in the case of piecewise $C^1$ continuous functions, equivalently, 
\begin{equation*}\label{hor-c-es-bis}
m(x,y)M(x,y) >0 ,
\end{equation*} 
then $\mF=\{f_1, f_2\}$ satisfies the nonsmooth H\"ormander condition of step $2$ at $(x,y,z)^t\in \R^3$. 

Now we can 
apply Theorem~\ref{eik-thm} 
 to the eikonal Dirichlet problem 
\begin{equation}\label{DP-ex-2}
\begin{cases} |\partial_x u +\a(y) \partial_z u|^2 + |\partial_y u + \beta(x) \partial_z u|^2 =1 & \mbox{in $\Omega$} ,\\
                     u =g & \mbox{on $\partial \Omega$ }  . \end{cases}                 
\end{equation}
{ We begin with the case $g\equiv 0$. If 
\eqref{hor-c-es} 
 holds 
at any 
 point $p =(x,y,z)^t$ 
 in $\Omega$ or 
 in $I(\partial \Omega)$, and  
 $\mT =\R^3 \setminus \Omega$ satisfies, for all $p \in \pT \setminus I(\pT)$, a $({\bf n}(p), \theta(p))$ inner cone condition  in a relative neighborhood of $p$ in $\pT$ 
with either 
$\mC_{\theta(p)}({\bf n}(p) ) \cap span\{f_1(p), f_2(p)\} \neq \{0\}$  
or $ \big(D_C\beta(x) - D_C \a(y)  \big) \partial_z \subset \mC_{\theta(p)}( {\bf n(p)})$, 
then  \eqref{DP-ex-2} 
admits a unique locally $(1/2)$-H\"older continuous viscosity solution on $\ol{\Omega}$.  All these conditions can be simplified as above if $\alpha, \beta$ are piecewise $C^1$ or $\partial\Omega \setminus I(\partial \Omega)$ is of class $C^1$.
  
When $g$ is a general continuous function on $\partial \Omega$, 
 the compatibility condition~\eqref{comp-con} is satisfied, for instance, 
  if $g$ is $C^{-1}$-Lipschitz on $\partial \Omega$ 
  with 
  $C = 1+ \sup\{ |\alpha(y)|  \vee  |\beta(x)|  :  (x,y, z)\in \ol{\Omega} \}$. }
  \end{example}

\begin{example}[Nonsmooth Gru\v sin type vector fields]\label{Grus-ex}
In $\R^n \times \R$ 
 consider the control system $\mF$ consisting of the vector fields 
\[
{ f_j = \partial_{x_j} \quad \mbox{ and } \quad f_{n+i} = \a_i(x_i) \partial_y, \quad \mbox{for } j=1, \ldots n , \;  i=1, \ldots m,}
\]
where $x= (x_1, \ldots, x_n) \in \R^n$, $(x, y) \in \R^n\times\R$, and  $\a_i : \R \to \R$ are  functions of class $C^{k-1, 1}$ for some $k\in \N$.  
We study the associated eikonal Dirichlet problem 
\begin{equation}\label{DP-ex}
\begin{cases} \sum_{i=1}^n|\partial_{x_i} u|^2 + 
\sum_{i=1}^m\a_i(x_i)^2 
 |\partial_y u|^2 =1 & \mbox{in $\Omega$} , \\
                     u =g & \mbox{on $\partial \Omega$}  , \end{cases}
\end{equation}
where $\Omega\subseteq \R^n \times \R$ is open  and $g$ is continuous.   
Albano \cite{Al12} studied this problem for { the classical Gru\v sin vector fields, where} $\alpha_i(x_i) =x_i^k$, and $g\equiv 0$ on the domain
\bel{dom-Alb}
\Omega = \{ (x, y) \in \R^n\times \R\; : \: y> M|x|^{k+1} \}, 
\eeq
for some $M>0$. 
He proved that the 
 unique viscosity  solution of the problem is locally $1/(k+1)$-H\"older continuous { in $\bar\Omega$}. 
 We  
extend  that result by considering more general $\Omega$, 
  $g$, and $ \a_i$ satisfying
 \[
 \a_i(x_i)=0 \quad \Longleftrightarrow \quad x_i=0 .
 \] 
 Then the vector fields {are of Gru\v sin type in the sense that} they span the whole space $\R^{n+1}$ at all points but those of the $y$-axis, $p=(0, \ldots, 0,y)$. At such points we need a non-null Lie bracket. 

We assume the following conditions.
 
 \noindent (i) $\partial\Omega \setminus I(\partial \Omega)$ is a $C^1$ manifold or the empty set; 
 
 \noindent  
{(ii)  
for any 
$p=(0, \ldots, 0,y)$ that either belongs to $\Omega\cup I(\partial \Omega)$, or belongs to $\partial\Omega \setminus I(\partial \Omega)$ and the outer normal ${\bf n}(p)$ to $\Omega$ at $p$ is parallel to the $y$-axis, for some $i$ 
either there exist  $j\leq k-1$ such that  
$
D^{j} \a_i(0) \ne 0
$
or 
\[
0 \notin D_C D^{k-1} \a_i(0) \, ;
\]
}

\noindent
 (iii) $g(p) -g(q) \le C^{-1}|p-q|$ or $g(p)-g(q) \le L^{-1}\ln{\left(1+ L|p-q|/C \right)}$ for all $p, q\in \partial \Omega$, where 
{ $C= 1+  \sup\{ |\a_i(x_i)| :  (x_1, \ldots, x_n, y) \in \overline{\Omega}, \; i=1, \ldots, n\}
$, and $L$ 
 is the maximum over  $i=1, \ldots, n$ of the  
 Lipschitz constants of the functions $\a_i$ on 
  the $i$-th projection of $\ol{\Omega}$. }   

Then problem \eqref{DP-ex} admits a unique continuous viscosity solution { bounded from below} 
which is in addition 
locally $(1/(k+1))$-H\"older continuous on $\overline{\Omega}$.

Indeed, 
this is a consequence of  Theorem~\ref{eik-thm}. 
At differentiability points of $\a_i$ one computes 
\[
\underbrace{[f_i, [f_i, [ \cdots [f_i,  f_{n+i} ]]]]}_{\text{$k$ bracketings }}(p) = D^{k}{ \a_i}(x_i) \partial_y , 
\] 
where $p=(x_1, \ldots, x_n, y)\in \R^n\times \R$. 
From this we deduce that  
\[
\underbrace{[f_i, [f_i, [ \cdots [f_i,  f_{n+i} ]]]]_{set}}_{\text{$k$ bracketings }}(p) = D_C D^{k-1}{ \a_i}(x_i) \partial_y , 
\]
Thus conditions (i) and (ii) guarantee that $\mF$ satisfies the nonsmooth H\"ormander's condition of step $k+1$ at any point of $\Omega \cup I(\partial \Omega)$, and 
$\mF$ is H-noncharacteristic of degree $k+1$ at any point of $\partial \Omega \setminus I(\partial \Omega)$. Condition (iii) guarantees that 
$g$ satisfy the compatibility condition ~\eqref{comp-con} via \eqref{comp-suff-1} or \eqref{comp-suff-2}.  
 Therefore 
  Theorem~\ref{eik-thm} applies to this example.
A simple explicit example of a nature not considered in the previous literature is the following:  $\Omega$ given by  \eqref{dom-Alb},
\begin{equation*}
\a_i(x_i) = \begin{cases} 2x_i^k/(1+ |x_i|^k) & \mbox{if }  x_i \ge0   \\
                                   x_i^k/(1+ |x_i|^k)  & \mbox{if }  x_i <0  \end{cases}  \quad \mbox{for all }  i=1, \ldots, n \,.  
\end{equation*}
and $g(p) =c |p|$ for $p\in \partial \Omega$ with $0\le c \le 1/3$.  Clearly conditions (i), (ii), (iii) are satisfied, and problem \eqref{DP-ex} has a unique continuous viscosity solution on $\Omega$ which is 
$1/(k+1)-$H\"older continuous on $\overline{\Omega}$.    
 \end{example}

\section{Necessary conditions for the H\"older continuity of the minimum time function} 
\label{necessary}

{ 
In this section we work in a different setting with respect to the rest of the paper: the system $\mF$ is not necessarily symmetric, but on the other hand we assume more regularity  on the target and the vector fields.

It is well-known that the Lipschitz continuity of the minimum time function is characterised by the existence of  vector fields pointing inward the target. More precisely, if $\partial \mT$ is $C^2$ near $x_0$, then 
\begin{equation}
\label{T-Hol-as-second}
T(x) \le C d(x)^\a ,
\end{equation}
 near $x_0$, with $\a>1/2$, if and only if there is $\bar f\in \mF$ such that
  { $\bar f(x_0) \cdot {\bf n}(x_0) < 0$, where ${\bf n}(x_0)$ is the outer normal unit vector of $\mT$ at $x_0$},
  see e.g. \cite[Theorem~5.5]{BaFa90}, { \cite{Vel}}. We also prove an extension of such result to set-valued systems of class $C^{-1,1}$ in the Appendix \ref{append}.
 
  In the main result of the section we want to characterise  the $\alpha$-H\"older continuity of $T$ 
in the range $\a \in(1/3,1/2]$. Hence we must restrict to the case $f(x_0) \cdot {\bf n}(x_0) \geq 0$ for all $f$.
We will assume the 
 stronger 
 condition:
 there exists $\nu>0$ such that for all $f\in \mF$
\begin{equation}
\label{tec-ass-2}
\mbox{either} \quad f(x_0) \cdot {\bf n}(x_0) \ge \nu \quad \mbox{or}\quad f(x_0) \cdot {\bf n}(x_0)= 0\,,
\end{equation}
which is automatically satisfied if $\mF$ is a finite set.

We know from Section \ref{sufficient} that a sufficient condition for the estimate \eqref{T-Hol-as-second} with $\a=1/2$ { in the case of symmetric systems} is the existence of a bracket 
pointing inward the target. Our main result states that such condition is also necessary if completed with the possibility of entering $\mT$ using a single tangential vector field.
}
{ Moreover, when we extend it to general systems, we also relax \eqref{T-Hol-as-second} to the weaker
\begin{equation}
\label{T-Hol-as}
T(x) \le C |x-x_o|^\a .
\end{equation}
}

First we need to recall a lemma on Taylor expansions of piecewise smooth trajectories. 
\begin{lem}\label{ex-pr-fl-le}
If $f_1, \ldots, f_m$ are $C^2$ vector fields on $\R^n$, then, for $s_1, \ldots, s_m \ge 0$,  
\begin{equation}\label{Tex-pr-fl-O}
\psi_{s_1}^{f_1} \circ \cdots \circ \psi_{s_m}^{f_m}(x) = x +  r(x)  +  R(x) + O(t^3)   
\end{equation}
as $t = \sum_{i=1}^m s_i\to 0$, where
\begin{equation}\label{def-r-R}
r(x) = \sum_{i=0}^m f_i(x)s_i 
\,, \quad R(x) =  \frac{1}{2}\Big( Dr(x) r(x) +  \sum_{\substack{ i, j =1 \\ i<j }}^m [f_i, f_j](x)s_i s_j \Big) \,,
\end{equation}
\end{lem}
The proof of the lemma relies on applying recursively Taylor expansions of the flows; more details can be found, e.g., in~\cite{li}.

Let $\mT\subset \R^n$ 
 be { the closure of an open} set with a nonempty boundary. 
If the boundary of $\mT$ is a $C^k$ manifold for $k\ge 2$, it is known,
see e.g. Gilbarg and Trudinger \cite[Section~14.6.]{GiTr01},
that the 
distance function $d(x) = {\rm dist}(x, {\partial\mT})$ extends on a $\delta$-neighborhood of $\pT$, $\d>0$, to a function of class $C^k$, which we denote still by $d$, such that $\nabla d(x_0) = {\bf n}(x_0)$ for all $x_0\in \pT$. 

\begin{thm}\label{nec-cond-2}
Let $\mT$ be the closure of an open set with $C^3$ boundary and $x_0\in  \pT$.  
Let $\mF$ be a 
control system 
uniformly of class $C^2$ 
in a neighborhood of $x_0$ { and satisfying \eqref{tec-ass-2}.
If for some $C>0$ and $\a\in ]1/3, 1/2]$ 
the estimate \eqref{T-Hol-as} holds
for all $x$ in a neighborhood of $x_0$,  then  either 
\begin{equation}
\label{trans-br1}
\exists \; f_1, f_2 \in \mF \; \text{such that}\quad
[f_1,f_2]\cdot {\bf n}(x_0)<0,
\end{equation}
or
\begin{equation}\label{trans-br2}
\exists \; \bar f\in \text{\upshape co} \{f\in\mF : f(x_0)\cdot {\bf n}(x_0)=0\} \; \text{such that}\quad\nabla(\nabla d\cdot \bar f)\cdot \bar f(x_0)<0 .
\end{equation}
}
\end{thm}
{ A consequence of the previous Theorem is a necessary and sufficient condition for symmetric systems. The fact that the implication in Theorem \ref{nec-cond-2} can be appropriately reversed also for a class of affine systems is proved in the paper of one of the authors \cite{so1}.}
{ \begin{corol}
Assume in addition that the set $\mF_o:= 
 \{f\in\mF : f(x_0)\cdot {\bf n}(x_0)=0\}$ is convex { and symmetric}. Then $T(x) \le C d(x)^{1/2}$ for all $x\notin {\mT}$ in a neighborhood of $x_0$ if and only if either \eqref{trans-br1} holds or  there is $\bar f\in \mF_o$ such that $\nabla(\nabla d\cdot \bar f)\cdot \bar f(x_0)<0$.
 \end{corol}
{\bf Proof.} The necessity part comes from the theorem above. The sufficiency of \eqref{trans-br1} for the H\"older continuity of $T$ follows from Theorem \ref{Hold}.
On the other hand, for any $f\in\mF_o$, a Taylor expansion in a neighborhood of $x_0$ gives
$$d(\psi^f_t(x_0))=\nabla(\nabla d\cdot f)\cdot f(x_0)\frac{t^2}2+O(t^3), \quad\text{as } t\to0^+ ,
$$
 and then the trajectory $\psi_.^{\bar f}(x_0) 
 $ associated  to the field $\bar f$ 
 enters the target in a time of order $t^2$. Standard continuous dependence with respect to the initial condition then shows that small time local attainability of the target holds at $x_0$ by means of the single vector field $\bar f$.
\qed}

\smallskip
{\bf Proof of Theorem \ref{nec-cond-2}.}
We assume by contradiction that,  for all $f_1, f_2\in \mF$, \eqref{trans-br1} and \eqref{trans-br2} do not hold, that is,
\begin{equation}\label{eqnec}
\begin{array}{ll}
[f_1, f_2](x_0) \cdot {\bf n}(x_0) { { =} 0 } ,\quad \mbox{for all } f_1, f_2\in \mF \mbox{ and }\\
{ \nabla(\nabla d\cdot g)\cdot g(x_0)\geq 0 , 
\quad\mbox{for all }g\in\mbox{co}{\mathcal F}_o .}
\end{array}\end{equation} 

Let $x_k = x_0 + {\bf n}(x_0) /k$, then for { sufficiently large $k$ we have that $|x_k -x_0| = d(x_k)$}. We choose a decreasing sequence $\e_k \to 0$ such that $  
0<\e_k \leq  T(x_k)^2 $ as $k \to \infty$.  For each $k$ we can find 
vector fields $f^k_1, \ldots, f^k_{m_k}$, nonnegative times $t_k$, $s^k_1, \ldots, s^k_{m_k}$ for some 
$m_k \in \N$ such that 
\[
 y_k=  \psi_{s^k_1}^{f^k_1} \circ \cdots \circ \psi_{s^k_{m_k}}^{f^k_{m_k} }(x_k) \in {\partial \mT}, 
\]
$\sum_{i=1}^{m_k}s_i^k = t_k $, and $ t_k < T(x_k) + \e_k$. 

From now on we drop the index $k$ from $x_k, t_k, f^k_i, s_i^k, m_k, y_k$.  By the previous lemma, the uniform $C^2$ regularity of the vector fields in $\mF$ around $x_0$,
and { $|x-x_0| = d(x)$}, we obtain 
\begin{equation}
\label{exp-y}
y= x+ r(x_0) + R(x_0) + O(t^3) + { O(t |x-x_0|)} \, 
\end{equation}
as $k\to \infty$ (hence, $t\to 0$, $d(x) \to 0$).
{ Note that $r(x_0)$ and $R(x_0)$  defined by \eqref{def-r-R} depend on $s_i^k$ and satisfy
\bel{rR}
r(x_0) = O(t)\,, \quad  Dr(x_0) r(x_0) = O(t^2)\,,  \quad R(x_0)= O(t^2)\,  \quad \text{as }\, t\to 0+.
\eeq
}
{By the first condition in \eqref{eqnec} }we get 
\begin{equation}
\label{exp-(x-y)n}
(x-y) \cdot {\bf n}(x_0) \le  -r(x_0) \cdot {\bf n}(x_0)  - { \frac 12 Dr(x_0) r(x_0) \cdot {\bf n}(x_0)} +   O(t^3) + { O(t |x-x_0| )}\, .
\end{equation}  
We now expand the distance function $d(\cdot)\in C^3$ around $x_0$ and obtain
\begin{equation}
\label{exp-d-1}
d(x) = {\bf n}(x_0) \cdot (x-x_0) + { O(|x-x_0|^2)} \, , 
\end{equation}
\begin{equation*}
{ 0  = }  
 d(y) = {\bf n}(x_0) \cdot (y-x_0) + \frac 12 D^2d(x_0)(y-x_0) \cdot (y-x_0) + O(|y-x_0|^3)  \, , 
\end{equation*}
From the last estimate, using the fact that 
$y- x =O(t)$ as $k \to \infty$, we deduce 
\begin{equation}
\label{exp-d-2}
{ 0 = }  
 d(y) = {\bf n}(x_0) \cdot (y-x_0) + \frac 12 D^2d(x_0)(y-x) \cdot (y-x) + O(t^3) + O(t |x-x_0|) + { O(|x-x_0|^2)} \, . 
\end{equation}
We subtract \eqref{exp-d-2} from \eqref{exp-d-1}, { then use (\ref{exp-(x-y)n}) with $r(x_0) \cdot {\bf n}(x_0) \ge 0$ 
{ to get 
\begin{align}
d(x) &\le {\bf n}(x_0) \cdot (x-y) -   \frac 12 D^2d(x_0)(y-x) \cdot (y-x) + O(t^3) + { O(t |x-x_0|) + O(|x-x_0|^2)}  \nonumber \\
&\le -\frac 12 Dr(x_0)r(x_0) \cdot {\bf n}(x_0)  - \frac 12 D^2d(x_0)(y-x) \cdot (y-x) 
 +  O(t^3) + { O(t |x-x_0|) + O(|x-x_0|^2)} . \nonumber 
\end{align}
Next we use (\ref{exp-y}), \eqref{rR}, 
 and the identity $\nabla (  \nabla d \cdot r )\cdot r = (Dr) r \cdot \nabla d +( D^2d)  r \cdot   r$
to obtain
\bel{dx-est}
{ |x-x_0|=}d(x) \leq - \frac 12 \nabla (  \nabla d \cdot r )(x_0) \cdot r( x_0) +   O(t^3) + { O(t  |x-x_0|) + O(|x-x_0|^2)} \quad\text{as $k \to \infty$} .
\eeq 
}
Also from \eqref{exp-(x-y)n}, { \eqref{rR}, and from} \eqref{exp-d-2} we obtain, respectively, 
\[
\begin{gathered}
r(x_0)\cdot {\bf n}(x_0) \leq (y-x)\cdot {\bf n}(x_0) +O(t^2) + { O(t|x-x_0|)}, \\
{ 0=} {\bf n}(x_0) \cdot (y-x_0) +  O(t^2) + { O(t |x-x_0|) + O(|x-x_0|^2)}\, ,
\end{gathered}
\]
and subtracting them we obtain  
\bel{exp-rn}
r(x_0)\cdot {\bf n}(x_0) \le O{ \left((t +{ |x-x_0|})^2\right) }\,.
\eeq
Now we introduce two sets of indices: $P_1=\{ i=1,\ldots, m\, : \, f_i(x_0) \cdot {\bf n}(x_0) =0\}$, $P_2 = \{1, \ldots, m\} \setminus P_1$, 
and split 
 $r(x_0) = r_1(x_0) + r_2(x_0)$, where $r_j(x_0) =\sum_{i\in P_j} f_i(x_0) s_i$, for $j=1,2$.
Note that when $P_1$ is not empty, { $r_1(x)=t_1 g(x)$, where $t_1=\sum_{i\in P_1}s_i$ and $g\in\mbox{co}\{f_j:j\in P_1\}\subseteq \mbox{co}{\mathcal F}_o$.}
By \eqref{tec-ass-2} 
\[
r(x_0)\cdot {\bf n}(x_0) = r_2(x_0)\cdot {\bf n}(x_0) \ge \nu \sum_{i\in P_2} s_i\,,
\]
and thus, by \eqref{exp-rn},
$\sum_{i\in P_2} s_i \le { O(t +|x-x_0|)^2}$, { which implies
\[
\nabla (  \nabla d \cdot r )(x_0) \cdot r( x_0)= \nabla (  \nabla d \cdot r_1 )(x_0) \cdot r_1( x_0) + { tO(t +|x-x_0|)^2+O((t+|x-x_0|)^4)} .
\]
Plugging this in \eqref{dx-est} and then using \eqref{eqnec} we get
\begin{equation*}\begin{array}{ll}
{ |x-x_0|} &\le   -\frac {t_1^2}2 \nabla (  \nabla d \cdot g )(x_0) \cdot g( x_0) +   { O(t^3) + O(t |x-x_0|) + O(|x-x_0|^2)}\\
&\leq { O(t^3 
+t |x-x_0| 
+|x-x_0|^2)\,. }
\end{array}\end{equation*}
}
Assumption \eqref{T-Hol-as} { gives { $|x-x_0|^{-1}=O(T(x)^{-1/\a})\leq O(t^{-1/\a})$}, and after dividing the previous inequality by { $|x-x_0|$} } we conclude
\[
1 \le O(t^{3 -1/\a} +t + { |x-x_0|} ) \quad\text{ as $k \to \infty$} ,
\]
which is a contradiction.  \qed 
\begin{rem}
The result extends to less regular control systems $\mF$ and targets $\mT$, and precisely, to $\mF$ of class $C^1$ and $\mT$ of class $C^2$ provided that the H\"older continuity assumption~\eqref{T-Hol-as} holds with $\alpha = 1/2$, instead of $\alpha >1/3$. 
The proof is similar to the one above, 
by using the less precise expansion 
\[
\psi_{s_1}^{f_1} \circ \cdots \circ \psi_{s_m}^{f_m}(x_0) = x_0 +  r(x_0)  +  R(x_0)  + o(t^2)  , \quad t = \sum_{i=1}^m s_i\to 0 ,
\]
 instead of \eqref{Tex-pr-fl-O}.   
\end{rem}

\subsection{Appendix: 
On Lipschitz continuity for $C^{-1,1}$ vector fields}
\label{append}
Here we extend to systems of class $C^{-1, 1}$ (that is, consisting of upper semi continuous set-valued vector fields with compact, convex, nonempty values) the necessary condition for local Lipschitz continuity of $T$ well-known for locally Lipschitz systems,
see~\cite[Theorem~5.5]{BaFa90}.  

\begin{lem}\label{rough-exp}
Let $f_1, \ldots, f_m$ be vector fields of class $C^{-1, 1}$ on $\R^n$ 
and $x_* \in \R^n$. 
Then for $s_1,   \ldots, s_m \ge 0$,  $x \in \R^n$ and 
\[
y\in \psi^{f_1}_{s_1} \circ \cdots \circ \psi^{f_m}_{s_m} (x)
\]
we have 
\[
{\rm dist }\, \left( y-x,\,  \sum_{i=1}^m f_i(x_*)s_i \right) = t \,\gamma ( O(t + |x-x_*| )   )   
\]
as $t + |x-x_*| \to 0$, where $t =\sum_{i=1}^m s_i$, and $\gamma$ is the sum of the upper semicontinuity moduli of the 
 vector fields $f_i$ 
 defined by \eqref{usc-mod}. 
\end{lem}
\noindent{\bf Proof.}  
 It follows by induction on 
  $m$, by applying repeatedly estimate \eqref{B-exp-m1},{ that we prove next.
Let $\rh, M >0$ be such that $f(x) \subset M B_1(0)$
for $|x-x_*| \le \rh$, and $\psi^f(x, t)$ is not empty for $|x-x_*|\le  \rh$, $|t| \le \rh$, see \cite{AuCe84}.  
Set $\rb = \rh/(M+2)$. Now take $(x, t)$ so that 
$\r= |t| +|x-x_*|\ \le \rb$.  For any $t\in \R$, we denote by $I_t$ the interval $[\min\{0,t\},\;\max\{0,t\}  ]$. 
Let $y \in \psi^f(x, t)$. There exists a trajectory $\xi : I_t  \to \R^n$ such that $\xi(0) =x $, $\xi(t) = y$, and ${\dot \xi}(s) \in f(\xi(s) )$
for almost every $s \in I_t$.
We show that $|\xi (s)  -x_*| \le \rh$ for all $s \in I_t$. If this were not the case, then there would exist a  $\tau \in I_t$
such that   $|\xi (s)  -x_*| \le \rh$ for all $s \in I_\tau$, and  $|\xi (\tau)  -x_*| = \rh$. 
Then 
\[
\rh= |\xi(\tau )-x_*| \le |\xi(\tau ) -x| + |x-x_*| \le \Big| \int_0^\tau {\dot \xi}(s)  ds \Big| + \r \le M|\tau| +\r \le (M+1)\r <\rh,   
\]
Therefore, it must be $\tau = t$, and we have also shown that 
\[
|\xi(s )-x_*| \le (M+1) (|s| + |x-x_*| ) \quad \mbox{whenever } s \in I_t.
\]
Thus we have 
\[
| y-x- tv| = \Big|\int_0^t ( {\dot \xi}(s)  -v ) ds \Big| \le |t| \gamma\big( (M+1) (|t| +|x-x_*|\ )   \big) 
\]
for all $v\in f(x_*)$, 
proving the desired estimate.}
 \qed  
\begin{thm}
Consider a control system $\mF$ of class $C^{-1,1}$ and 
 $\mT\subset \R^n$ 
  of class $C^1$ around a point $x_0 \in \pT$. If for some $C\ge 0$ the estimate \eqref{T-Hol-as} with $\a=1$ holds
in a neighborhood of $x_0$, then there exist $f \in \mF$ and $v\in f(x_0)$ such that 
 \bel{out-nor}
  {\bf n }(x_0) \cdot v <0  \, ,
 \eeq
 where ${\bf n }(x_0)$ is the outer normal of $\mT$ at $x_0$. 
\end{thm}
\noindent{\bf Proof.}  Assume by contradiction that ${\bf n }(x_0) \cdot v \ge 0$ for all $f\in \mF$, $v\in f(x_0)$. 
Take $x_k = x_0+  {\bf n }(x_0)/k$ and 
a sequence of times $t_k$ such that $t_k <T(x_k) + \e_k$, for some decreasing sequence $\e_k\to 0$ such that $\e_k =o(T(x_k))$ as $k\to \infty$.  
{ Then the assumption $T(x)\leq Cd(x)$ implies
\bel{nuova}
t_k =O(d(x_k) ) .
\eeq
}
By Lemma~\ref{rough-exp} and the definition of $T(x_k)$ we can find points $y_k \in \mT$ such that 
\bel{exp-yk}
y_k = x_k + \sum_{i=1}^{m_k} v_i^ks_i^k + o(t_k + |x_k-x_0| ) 
\eeq
as $k \to \infty$, for suitable $f_i^k \in \mF$, $v_i^k \in f_i^k(x_0)$, $s_i^k \ge 0$ such that $\sum_{i=1}^{m_k}{s_i^k} =t_k$. 
Since ${\bf n }(x_0) \cdot f(x_0) \subset \R_+$ for all $f\in \mF$ and $|x_k{-x_0}| \le Cd(x_k)$ for all $k\in\N$, we obtain 
\bel{exp-ykxk}
(x_k - y_k) \cdot {\bf n}(x_0) \le  o(t_k +d(x_k))\, .
\eeq   

As in the proof of Theorem~\ref{Hold}, we can find a function $w$ which is differentiable at $x_0$ and a constant $c\ge 0$ such that   $d(x) \le w(x) \le c\, d(x)$ for all $x$ in a neighborhood of $x_0$, $x\notin \mT$,  and  
$\nabla w(x_0) = {\bf n}(x_0)$, $w(y) \le 0$ for all $y$ in a neighborhood of $x_0$, $y\in \mT$.  Expanding that function at $x_0$ we find 
\[
\begin{gathered}
w(x_k) = {\bf n}(x_0) \cdot (x_k-x_0) + o(|x_k-x_0|)\,,\\ 
0\ge  w(y_k) = {\bf n}(x_0) \cdot (y_k-x_0) + o(|y_k-x_0|)\,.
\end{gathered}
\]
Subtracting the two expansions, using the estimates $|x_k{-x_0}| \le Cd(x_k)$, 
$|y_k - x_k| \le C t_k + o( d(x_k)) $  following from expansion \eqref{exp-yk}, and 
\eqref{exp-ykxk}, we obtain
\[
d(x_k) \le o(t_k + d(x_k)  )
\]
as $k\to \infty$. 
{ This is  a contradiction to \eqref{nuova}}. \qed 

\bigskip
\noindent
{
{\it Acknowledgements.} We wish to thank Franco Rampazzo for useful conversations and both referees for insightful comments.
}
\smallskip

\bibliographystyle{plain}

\end{document}